\documentclass[11pt,twoside]{article}
\usepackage{amsthm,amsfonts,mathtools,amscd,amssymb, mathrsfs}
\usepackage{euscript}
\usepackage{enumitem}
\usepackage{graphicx}
\usepackage{tikz}
\usepackage[font=small]{caption}
\usepackage{subcaption}

\usepackage{cite}
\usepackage[utf8]{inputenc}
\usepackage[english]{babel}
\usepackage{xcolor}

\usepackage{jmpag-n} 

\hypersetup{linkcolor=darkblue, urlcolor=darkred, citecolor=darkgreen}

\begin{document}


\title[Controllability Problems for the Heat Equation on a Half-Plane]%
{Controllability Problems for the Heat Equation on a Half-Plane Controlled by the Neumann Boundary Condition with a Point-Wise Control}

\author{Larissa  Fardigola}
\address{B. Verkin Institute for Low Temperature Physics and Engineering of the
	National Academy of Sciences of Ukraine, 47 Nauky Ave., Kharkiv, 61103, Ukraine,\\
	V.N. Karazin Kharkiv National University, 4 Svobody Sq., Kharkiv, 61022, Ukraine}
\email{fardigola@ilt.kharkov.ua}

\author{Kateryna  Khalina}
\address{B. Verkin Institute for Low Temperature Physics and Engineering of the
	National Academy of Sciences of Ukraine, 47 Nauky Ave., Kharkiv, 61103, Ukraine}
\email{khalina@ilt.kharkov.ua}


\BeginPaper 


\newcommand{\LL}{L^2(\mathbb R_+)}

\newcommand{\HH}{\EuScript{H}}

\newcommand{\RR}{\mathbf{R}_T}

\newcommand{\RRL}{\mathbf{R}_T^L}
\newcommand{\fnl}{\left[\kern-0.1em\left]}
\newcommand{\fnr}{\right[\kern-0.1em\right]}
\newcommand{\sgn}{\mathop{\mathrm{sgn}}}
\newcommand{\supp}{\mathop{\mathrm{supp}}}

\newcommand{\zco}{{\text{$\bigcirc$\kern-6.2pt\raisebox{-0.5pt}{$0$}}}}
\newcommand{\zcx}[1]{~\mbox{#1\kern-7.5pt\raisebox{0.5pt}{$\bigcirc$}}}



\begin{abstract}
In the paper, the problems of controllability and approximate
controllability are 
studied for the control system $w_t=\Delta w$, $w_{x_1}(0,x_2,t)=u(t)\delta(x_2)$,
$x_1>0$, $x_2\in\mathbb R$, $t\in(0,T)$, where $u\in L^\infty(0,T)$ is a control. 
To this aid, it is investigated the set $\EuScript R_T(0)\subset L^2((0,+\infty)\times\mathbb R)$ of its end states which are reachable from 0. It is established that a function  $f\in\EuScript R_T(0)$ can be represented in the form $f(x)=g\big(|x|^2\big)$ a.e. in $(0,+\infty)\times\mathbb R$ where $g\in L^2(0,+\infty)$.  In fact, we reduce the problem dealing with functions from $L^2((0,+\infty)\times\mathbb R)$ to a problem dealing with functions from $L^2(0,+\infty)$.
Both a necessary and
sufficient condition for  controllability and a sufficient condition for approximate controllability in a given time $T$ under a control $u$ bounded by a given constant are
obtained in terms of solvability of a Markov power moment problem.  
Using the  Laguerre functions (forming an orthonormal basis of $L^2(0,+\infty)$), necessary and sufficient conditions for approximate controllability and numerical solutions to the approximate controllability problem are obtained. It is also shown that there is no initial state that is null-controllable in a given time $T$. The results are illustrated by an example.

\key{heat equation, controllability, approximate controllability, half-plane}

\end{abstract}



\section{Introduction}\label{sect1}

In the paper,
the controllability problems for the heat equation are studied on a half-plane. Note that these problems for the heat equation were studied both in  bounded and unbounded domains. However, most of the papers studying these problems deal with domains bounded with respect to the spatial variables (see some recent papers: \cite{ZK, WL,  QL, DEM, BLM, B-P-1, O, GYZ}, and the references therein). At the same time, there are quite a few papers considering domains unbounded with respect to the spatial variables \cite{MenCab, CabMenZua, GonTer, Ter, TerZua, CMV, Mill, Fat, MZua1, MZua2, FKh, FKh2, Barb, DWZh, DD, WWZZ, FKh3, FKh4, FKh5, BGST}.

A point-wise control is a mathematical model of a source supported in a domain of very small size with respect to the while domain. That is why studying control problems under a point-wise control is an important issue in control theory (see, e.g. \cite{NN, CF, JLL, HS, DR} and others).

In \cite{MZua2}, the boundary controllability of the 2-d heat equation was studied in a half-space. Using similarity variables and weighted Sobolev spaces and developing solutions in Fourier series reduce the control problem to a sequence of one-dimensional
controlled systems. The null-controllability properties of these systems had been
studied  in \cite{MZua1}. It had been proved that no initial datum belonging to any Sobolev space of negative order may be driven to zero in finite time. In \cite{MZua2}, it was established  that if all the corresponding 1-d problems are null-controllable, then the multidimensional problem is null-controllable. However, it was also proved that if there exists at least one 1-d problem which is not null-controllable,
then the multi-dimensional problem is not null-controllable. The results of the one-dimensional case was applied to obtain the corresponding results for the multi-dimensional case.

The  controllability problems for the heat equation on a half-plane controlled by the Dirichlet boundary condition with  a point-wise control were studied in \cite{FKh3}. 
Both necessary
and sufficient conditions for  controllability and sufficient  conditions for approximate controllability in a given time  under a control  bounded by a
given constant were obtained in terms of solvability of a Markov power moment problem. Orthogonal bases in special spaces of Sobolev type (consisting of functions of two variables) were
constructed by using the generalized Laguerre polynomials. Applying these bases, necessary and sufficient conditions for approximate
controllability and numerical solutions to the approximate controllability
problem were obtained.

The boundary controllability of the wave equation on a half-plane $x_1>0$, $x_2\in \mathbb R$ with a pointwise control on the boundary  was studied in \cite{FLVJMPAG05, FLVJMPAG15, LVF0}.

Consider the following control system on a half-plane
\begin{align}
&w_t= \Delta w,&& x_1>0,\ x_2\in\mathbb{R},\ t\in(0,T),\label{eq}\\ 
&w_{x_1}\big(0,(\cdot)_{[2]},t\big)=\delta_{[2]}u(t),&& x_2\in\mathbb{R},\ t\in(0,T),\label{bc}\\
&w\big((\cdot)_{[1]},(\cdot)_{[2]},0\big)=w^0,&&x_1>0,\ x_2\in\mathbb{R}, \label{ic}
\end{align}
where $\Delta =(\partial/\partial x_1)^2+(\partial/\partial x_2)^2$, $T>0$, $u\in L^\infty(0,T)$ is a control, $\delta_{[m]}$ is the Dirac distribution with respect to $x_m$, $m=1,2$. The subscripts $[1]$ and $[2]$ associate with the variable numbers, e.g., $(\cdot)_{[1]}$ and $(\cdot)_{[2]}$ correspond to $x_1$ and $x_2$, respectively, if we consider $f(x)$, $x\in\mathbb R^2$. This control system is considered in spaces of Sobolev type (see details in Section \ref{sect2}). We treat equality \eqref{bc} as the value of the distribution $w_{x_1}$ on the line $x_1=0$ (see Definition \ref{defpoint} below).

In Section \ref{sect2}, some notation and definitions are given. 

In Section \ref{sect3}, control problem \eqref{eq}--\eqref{ic}  is reduced to control problem \eqref{eq1}, \eqref{ic1} (see below) by using the even extension with respect to $x_1$ for $w(\cdot,t)$ and $w^0$, $t\in[0,T]$. It is proved that systems \eqref{eq}--\eqref{ic} and \eqref{eq1}, \eqref{ic1} are equivalent so, basing on this reason, we consider control system \eqref{eq1}, \eqref{ic1} (dealing with functions defined on $\mathbb R^2$)  instead of control system \eqref{eq}--\eqref{ic} (dealing with functions defined on $[0,+\infty)\times\mathbb R$).
The set $\EuScript R_T(0)\subset L^2\big(\mathbb R^2\big)$ of its states reachable from 0 (i.e. the set which is formed by the end states $w(\cdot,T)$ of control system \eqref{eq1}, \eqref{ic1} when controls $u\in L^\infty(0,T)$) and the set $\EuScript R_T^L(0)\subset\EuScript R_T(0)\subset L^2\big(\mathbb R^2\big)$ of its states reachable from 0 by using the controls $u\in L^\infty(0,T)$ satisfying the restriction $\|u\|_{L^\infty(0,T)}\leq L$ (where $L>0$ is a given constant)  are studied. 
In particular, properties of the solutions (Theorem \ref{th-sol}) and properties of the reachability sets $\EuScript R_T(0)$ and $\EuScript R_T^L(0)$ (Theorem \ref{reachprop}) are proved for this system.  
It is also established that a function  $f\in\EuScript R_T(0)$ can be represented in the form $f(x)=g\big(|x|^2\big)$ a.e. in $\mathbb R^2$ where $g\in L^2(0,+\infty)$. Therefore, the functions $g$ form the dual sets $\RR$ and  $\RR^L$ for the sets $\EuScript R_T(0)$ and $\EuScript R_T^L(0)$, respectively. In fact, the problem dealing with functions from $L^2\big(\mathbb R^2\big)$ is reduced to a problem dealing with functions from $L^2(0,+\infty)$. To this aid,  operators $\Psi$ and  $\Phi$ are introduced and studied. The results mentioned above  are applied in Sections \ref{sect3}--\ref{sect6}.  
In Section \ref{sect3} the following assertions are formulated for system \eqref{eq1}, \eqref{ic1}:
\begin{enumerate}[label={\upshape \arabic*)}, leftmargin=4ex, labelwidth=4ex, itemsep=0pt]
	\item
	some additional properties of the set  $\EuScript R_T^L(0)$ (Theorems \ref{thnec}--\ref{thmom}, \ref{thmomap}, and \ref{thclos}); 
	\item  
	necessary and sufficient conditions for controllability in a given time under the control bounded by a given constant (Corollary \ref{thmom-cor});
	\item
	sufficient conditions for approximate controllability in a given time under the control bounded by a given constant  (Corollary \ref{thmomap-cor});
	\item
	necessary and sufficient conditions for approximate controllability in a given time  (Theorem \ref{thappc});
	\item
	the lack of controllability to the origin  (Theorem \ref{thnulcontr}).
\end{enumerate}

In Section \ref{prop}, properties of the sets $\RR$ and  $\RR^L$ are established (Theorems \ref{pthnec}, \ref{pentire}, \ref{pcond-f}--\ref{pthmomap}, and \ref{pthclos}). 
In the proof of Theorem \ref{pthclos}, an algorithm for construction of controls solving the approximate controllability problem for system \eqref{eq1}, \eqref{ic1} is given.

In Section \ref{ex}, Theorem \ref{thappc} is illustrated by an example.

The results of  Section  \ref{prop} are applied in the proofs of Theorems
\ref{thnec}--\ref{thmom}, \ref{thmomap}, and \ref{thclos}  in Section \ref{sect6}. In this section Theorems \ref{th-sol} and \ref{thnulcontr} are also proved.

The main results of the present paper are rather similar to those of \cite{FKh3}. However, the methods of obtaining them are essentially different in these two papers. Roughly speaking, we deal with the two-dimensional case  studying reachability sets and constructing the solutions to controllability and approximate controllability problems in \cite{FKh3}, but reducing the two-dimensional reachability sets to the one-dimensional ones, we deal with the one-dimensional case studying these problems and constructing their solutions in the present paper. In addition, the methods used to study the one-dimensional reachability sets in this paper principally differ from those used for two-dimensional sets in \cite{FKh3}. That is why  Theorems \ref{thmom}, \ref{thmomap} and Corollaries \ref{thmom-cor}, \ref{thmomap-cor} in the present paper  also  differ from their analogues from  \cite{FKh3}. Moreover, Theorems \ref{thnec}--\ref{cond-f} have not analogues in \cite{FKh3}.

\section{Notation and preliminary results}\label{sect2}

Let $n\in\mathbb N$. By $|\cdot|$, we denote the Euclidean norm in $\mathbb{R}^n$.

Let $\mathscr{S}(\mathbb R^n)$ be the Schwartz space of rapidly decreasing functions \cite{Schw}. Put $\mathscr{S}=\mathscr{S}(\mathbb{R})$. Let $\mathscr{S}'(\mathbb{R}^n)$ ($\mathscr{S}'$) be the dual space for $\mathscr{S}(\mathbb{R}^n)$ ($\mathscr{S}$, respectively). 
Denote $\mathbb{R}_+=(0,+\infty)$. Let $\mathscr D(\mathbb R_+)$ be the space of infinitely differentiable functions on $\mathbb R$ whose supports are compact and they are contained in $\mathbb R_+$.

Let $D=\big(-i\partial/\partial x_1,\ldots,-i\partial/\partial x_n\big)$, $D^\alpha=\big(-i(\partial/\partial x_1)^{\alpha_1},\ldots,-i(\partial/\partial x_n)^{\alpha_n}\big)$, where $\alpha=(\alpha_1,\ldots,\alpha_n)\in\mathbb{N}^n_0$ is multi-index, $|\alpha|=\alpha_1+\cdots+\alpha_n$, $\alpha!=\alpha_1!\cdots\alpha_n!$, $\mathbb{N}_0=\mathbb{N}\cup\{0\}$. 

Consider the following Sobolev spaces \cite[Chap.~1]{VG}
\begin{equation*}
H^s(\mathbb{R}^n)=\left\{\varphi\in L^2(\mathbb{R}^n)\mid \forall\alpha\in\mathbb{N}_0^n \quad \big(|\alpha|\leq s \Rightarrow D^\alpha\varphi\in L^2(\mathbb{R}^n)\big)\right\},\quad s=\overline{0,2},
\end{equation*}
with the norm
\begin{equation*}
\left\| \varphi\right\|^s=\left(\sum_{|\alpha|\leq s}
\frac{s!}{(s-|\alpha|)!\alpha!}
\left(\left\| D^\alpha\varphi\right\|_{L^2(\mathbb{R}^n)} \right)^2\right)^{1/2},\quad\varphi\in H^s(\mathbb{R}^n).
\end{equation*}
We have $H^{-s}(\mathbb{R}^n)=\left( H^s\right(\mathbb{R}^n))^*$ with the norm $\left\|\cdot\right\|^{-s}$ associated with the strong topology of the adjoint space. Evidently, $H^0(\mathbb{R}^n)=L^2(\mathbb{R}^n)=\left( H^0(\mathbb{R}^n)\right)^*$.

For $s=\overline{0,2}$, denote
\begin{align*}
H_\zco^s=&\left\{\varphi\in L^2(\mathbb{R}_+\times\mathbb{R})\ \middle|\  \Big( \forall\alpha\in\mathbb{N}_0^2\quad \big(\alpha_1+\alpha_2\leq s\Rightarrow D^\alpha\varphi\in L^2(\mathbb{R}_+\times\mathbb{R})\big)\Big)\right.\\
&\kern31ex\left.\wedge  \Big(\forall k=\overline{0,s-1} \quad  D^{(k,0)}\varphi(0^+,(\cdot)_{[2]})=0\Big)\right\}
\end{align*}
with the norm
\begin{equation*}
\left\| \varphi\right\|_\zco^s=\left(\sum_{\alpha_1+\alpha_2\leq s} 
\frac{s!}{(s-(\alpha_1+\alpha_2))!\alpha_1!\alpha_2!}
\left(\left\| D^\alpha\varphi\right\|_{L^2(\mathbb{R}_+\times\mathbb{R})} \right)^2\right)^{1/2},\quad \varphi\in H_\zco^s\,,
\end{equation*}
and $H_\zco^{-s}=\left( H_\zco^s\right)^*$ with the norm 
$\left\|\cdot\right\|_\zco^{-s}$ associated with the strong topology of the adjoint space. Obviously, $H_\zco^0=L^2(\mathbb{R}_+\times\mathbb{R})=\left( H_\zco^0\right)^*$.

Consider also the spaces \cite[Chap.~1]{VG}
\begin{equation*}
H_m(\mathbb{R}^n)=\left\{ \psi\in L_{\text{loc}}^2(\mathbb{R}^n)\ \middle|\  \left(1+|\sigma|^2 \right)^{m/2}\psi\in L^2(\mathbb{R}^n)\right\},\quad m=\overline{-2,2},
\end{equation*}
with the norm
\begin{equation*}
\left\| \psi\right\|_m=\left\|\left(1+|\sigma|^2 \right)^{m/2}\psi \right\|_{L^2(\mathbb{R}^n)},\quad\psi\in H_m(\mathbb{R}^n).
\end{equation*}
Evidently, $H_{-m}(\mathbb{R}^n)=\left( H_m(\mathbb{R}^n)\right)^*$. It is easy to see that $H_0(\mathbb{R}^n)=H^0(\mathbb{R}^n)$. 

By $\langle f,\varphi\rangle$, denote the value of a distribution $f\in \mathscr S'(\mathbb R^n)$ on a test function $\varphi \in \mathscr S(\mathbb R^n)$. 

Let $\mathscr F: \mathscr S'(\mathbb R^n)\to \mathscr S'(\mathbb{R}^n)$ be the Fourier transform operator
with the domain $\mathscr S'(\mathbb R^n)$. This operator is an extension of the classical Fourier transform operator and is given by the formula
\begin{equation*}
\langle \mathscr F f,\varphi\rangle=\langle f,\mathscr F^{-1}\varphi\rangle,\quad f\in \mathscr S'(\mathbb R^n),\ \varphi\in \mathscr S(\mathbb R^n).
\end{equation*}
Due to \cite[Chap.~1]{VG}, the operator $\mathscr F$ is an isometric isomorphism of $H^m(\mathbb{R}^n)$ and $H_m(\mathbb{R}^n)$, $m=\overline{-2,2}$.

In the spaces $H_m(\mathbb{R}^2)$ and $H^m(\mathbb{R}^2)$, $m=\overline{-2,2}$, we consider the following inner products 
\begin{align*}
&\langle f,g \rangle_m=\left\langle\left(1+|\sigma|^2\right)^{m/2}f,\left(1+|\sigma|^2\right)^{m/2} g\right\rangle_0,&& f\in H_m(\mathbb{R}^2),\ g\in H_m(\mathbb{R}^2),\\
&\langle h,p\rangle^m=\langle \mathscr Fh,\mathscr Fp\rangle_m,&& h\in H^m(\mathbb{R}^2),\ p\in H^m(\mathbb{R}^2),
\end{align*}
where $\langle\cdot,\cdot\rangle_0$ is the inner product in $L^2(\mathbb{R}^2)$. Note that $\langle\cdot,\cdot\rangle^0=\langle\cdot,\cdot\rangle_0$.

A distribution $f\in \mathscr S'(\mathbb R^2)$  is said  to be
\emph{odd with respect to} $x_1$, if $\big\langle f,\varphi\big((\cdot)_{[1]},(\cdot)_{[2]}\big)\big\rangle=-\big\langle f,\varphi\big(-(\cdot)_{[1]}, (\cdot)_{[2]}\big)\big\rangle$, where $\varphi\in \mathscr S(\mathbb R^2)$.
A distribution $f\in\mathscr S'(\mathbb R^2)$  is said  to be
\emph{even with respect to} $x_1$, if $\big\langle f,\varphi\big((\cdot)_{[1]},(\cdot)_{[2]}\big)\big\rangle=\big\langle f,\varphi\big(-(\cdot)_{[1]},(\cdot)_{[2]}\big)\big\rangle$, where $\varphi\in \mathscr S(\mathbb R^2)$. 

Let $m=\overline{-2,2}$. By $\widehat{H}^m\left(\mathbb{R}^2\right)$ (or $\widehat{H}_m\left(\mathbb{R}^2\right)$), denote the subspace of all distributions in
$H^m\left(\mathbb{R}^2\right)$ (or $H_m\left(\mathbb{R}^2\right)$, respectively) that are even with respect to $x_1$. Evidently, $\widehat{H}^m\left(\mathbb{R}^2\right)$ (or $\widehat{H}_m\left(\mathbb{R}^2\right)$) is a closed subspace of $H^m\left(\mathbb{R}^2\right)$ (or $H_m\left(\mathbb{R}^2\right)$, respectively).


Let $f\in L^2\big(\mathbb{R}^2\big)$. Let also $f(x)=g(|x|^2)$, $x\in\mathbb R^2$, where $g$ is a function defined on $\mathbb R_+$. Setting $x_1=\sqrt r \cos \phi$, $x_2=\sqrt r \sin \phi$, $r\in\mathbb R_+$, $\phi\in[0,2\pi)$, we get
\begin{align}
\label{normh}
\|f\|_{L^2 (\mathbb{R}^2 )}
&=\left(\iint_{\mathbb R^2} |f(x)|^2\, dx\right)^{1/2}
=\left(\iint_{\mathbb R^2} |g(|x|^2)|^2\, dx\right)^{1/2}
\notag\\
&=\left(\pi\int_0^\infty |g(r)|^2\, dr\right)^{1/2}
=\sqrt\pi \|g\|_{\LL}.
\end{align}

Thus, if $f\in L^2\big(\mathbb{R}^2\big)$ and $f(x)=g(|x|^2)$, $x\in\mathbb R^2$, for some $g$  defined on $\mathbb R_+$, then $g\in\LL$ and \eqref{normh} holds; and vice versa: if $g\in\LL$, then for $f(x)=g(|x|^2)$, $x\in\mathbb R^2$, we have $f\in L^2\big(\mathbb{R}^2\big)$.

Taking this into account, we can introduce the space 
\begin{equation}
\label{defh}
\EuScript{H}=\left\{f\in L^2\big(\mathbb{R}^2\big)\ \middle|\ \exists g\in L^2(\mathbb R_+)\quad f(x)=g(|x|^2)\ \text{a.e. on}\ \mathbb{R}^2\right\}
\end{equation}
and the operator $\Psi: \HH\to\LL$ with the domain $D\big(\Psi\big)=\HH$ for which
$$
\Psi f=g \Leftrightarrow \big(f(x)=g(|x|^2) \ \text{a.e. on}\ \mathbb{R}^2\big),\quad f\in D\big(\Psi\big)=\HH.
$$
One can see that $\Psi$ is invertible, $\Psi^{-1}: \LL\to\HH$, and $\big(\Psi^{-1}g\big)(x)=g(|x|^2)$, $x\in\mathbb R^2$ for $g\in D\big(\Psi^{-1}\big)=\LL$. 

Summarising, we obtain the following proposition. 

\begin{proposition}
	\label{prop-1}
	The following assertions hold:
	\begin{enumerate}[label={\upshape (\roman*)}, ref={\upshape\thetheorem\ (\roman*)}, leftmargin=5ex, labelwidth=5ex]
		\item \label{prop-1-i}
		$\Psi$ is an isomorphism of $\HH$ and $\LL$\textup{;}
		\item \label{prop-1-ii}
		$\HH$ is a subspace of $\widehat{H}^0\left(\mathbb{R}^2\right) =\widehat{H}_0\left(\mathbb{R}^2\right)\subset L^2\big(\mathbb{R}^2\big)$\textup{;}
		\item \label{prop-1-iii}
		$\HH$ is a Hilbert space with respect to the inner product $\langle\cdot, \cdot \rangle_{L^2(\mathbb{R}^2)}$\textup{;}
		\item \label{prop-1-iv}
		$\langle  f,  h \rangle_{L^2(\mathbb{R}^2)}
		=\pi \langle \Psi f, \Psi h \rangle_{\LL}$, $f\in\HH$, $h\in\HH$\textup{;}
		\item \label{prop-1-v}
		$\|\Psi\|=1/\sqrt \pi$\textup{;}
		\item \label{prop-1-vi}
		$\mathscr F \HH= \HH$.
	\end{enumerate}	
\end{proposition}

Let $f\in\HH$, $F=\mathscr F f$, $g=\Psi f$, and $G=\Psi F$. Then $G=\Psi \mathscr F \Psi^{-1} g$.

Let us introduce the operator $\Phi:\LL\to\LL$ with the domain $D(\Phi)=\LL$ by the rule
\begin{equation*}
\Phi g = \Psi \mathscr F \Psi^{-1} g,\quad g\in D(\Phi)=\LL.
\end{equation*}
Since the operator $\Psi$ is an  isomorphism of $\HH$ and $\LL$ (see  Proposition \ref{prop-1}) and $\mathscr F f= \mathscr F^{-1} f$ for $f\in\HH$, we conclude that $\Phi$ is invertible and $\Phi^{-1}=\Phi$, in particular, $\Phi$ is an isometric isomorphism of $\LL$.

Let us find a formula for calculating $\Phi g$ if $g\in\LL$. Put $G=\Phi g$. Setting 
$x_1=\sqrt r \cos \phi$, $x_2=\sqrt r \sin \phi$, $r\in\mathbb R_+$, $\phi\in[0,2\pi)$, and
$\sigma_1=\sqrt \rho \cos \theta$, $\sigma_2=\sqrt \rho \sin \theta$, $\rho\in\mathbb R_+$, $\theta\in[0,2\pi)$,
we get
\begin{align}
\label{fourier}
G(\rho)=\big(\mathscr F \Psi^{-1} g\big)(\sigma)
&=\frac1{2\pi} \lim_{N\to\infty}
\iint_{|x|\leq  N^2} \big(\Psi^{-1} g\big)(x) 
e^{-i\langle x,\sigma\rangle}\, dx
\notag\\
&=\frac1{4\pi} \lim_{N\to\infty} 
\int_0^N g(r)\int_0^{2\pi} e^{-i\sqrt{r\rho}\cos\phi} \,d\phi\,dr
\notag\\
&=\frac12 \lim_{N\to\infty} \int_0^N g(r) J_0\big(\sqrt{r\rho}\big)\, dr,\quad \rho\in \mathbb R_+,
\end{align}
where $J_0$ is the Bessel function of order $0$.
Here the relation 
$$
J_0(\xi)=\frac1{2\pi} \int_0^{2\pi} e^{-i\xi\cos\phi}\,d\phi, \quad \xi \in\mathbb R,
$$
has been used.

Summarising, we obtain the following proposition.

\begin{proposition}
	\label{prop-2}
	The following assertions hold:
	\begin{enumerate}[label={\upshape (\roman*)}, ref={\upshape\thetheorem\ (\roman*)}, leftmargin=5ex, labelwidth=5ex]
		\item \label{prop-2-i}
		$\Phi$ is invertible and $\Phi^{-1}=\Phi$\textup{;}
		\item \label{prop-2-ii}
		$\displaystyle \big(\Phi g\big)(\rho)=\frac12 \lim_{N\to\infty} \int_0^N g(r) J_0\big(\sqrt{r\rho}\big)\, dr$, $\rho\in\mathbb R_+$, $g\in\LL$.
	\end{enumerate}
\end{proposition}

With regard to Proposition \ref{prop-2-ii}, one can see that the transform providing by the operator $\Phi$ is a modification of the well-known Hankel transform of order~0.

Let $g\in H_{-s}(\mathbb R^2)$, $s=\overline{0,2}$. By a similar reasoning to that of \cite{FKh3}, we get
$g\big(\sigma_1, (\cdot)_{[2]}\big)\in H_{-s}(\mathbb R)$ for almost all $\sigma_1\in\mathbb R$, and  $g\big( (\cdot)_{[1]}, \sigma_2\big)\in H_{-s}(\mathbb R)$ for almost all $\sigma_2\in\mathbb R$.

Let $\psi\in H_s (\mathbb R)$, $s=\overline{0,2}$. Denote
\begin{align*}
\langle g,\psi \rangle_{[1]}(\sigma_2)
&=\int_{-\infty}^\infty \left(1+\sigma_1^2\right)^{-s}g(\sigma_1,\sigma_2)\left(1+\sigma_1^2\right)^s\overline{\psi(\sigma_1)}\,d\sigma_1,&& \sigma_2\in\mathbb R,
\\
\langle g,\psi \rangle_{[2]}(\sigma_1)
&=\int_{-\infty}^\infty \left(1+\sigma_2^2\right)^{-s}g(\sigma_1,\sigma_2)\left(1+\sigma_2^2\right)^s\overline{\psi(\sigma_2)}\,d\sigma_2, && \sigma_1\in\mathbb R.
\end{align*}
Then,
$$
\langle g,\psi \rangle_{[1]}\in H_{-s}(\mathbb R)\quad\text{and}\quad \langle g,\psi \rangle_{[2]}\in H_{-s}(\mathbb R).
$$

Let $f\in H^{-s}(\mathbb R^2)$, $\varphi\in H^s (\mathbb R)$, $s=\overline{0,2}$. Denote
$$
\langle f,\varphi \rangle_{[1]}
=\mathscr F_{\sigma_2\to x_2}^{-1}\left(\langle \mathscr F_{x\to\sigma}f,\mathscr F\varphi \rangle_{[1]} \right)\ \text{and}\ 
\langle f,\varphi \rangle_{[2]}
=\mathscr F_{\sigma_1\to x_1}^{-1}\left(\langle \mathscr F_{x\to\sigma}f,\mathscr F\varphi \rangle_{[2]} \right).
$$

Since the operator $\mathscr F$ is an isometric isomorphism of $H^m(\mathbb{R}^n)$ and $H_m(\mathbb{R}^n)$, $m=\overline{-2,2}$,  
we get
$$
\langle f,\varphi \rangle_{[1]}\in H^{-s}(\mathbb R)\quad\text{and}\quad \langle f,\varphi \rangle_{[2]}\in H^{-s}(\mathbb R).
$$

The following definition is given with regard to the definition of a distribution's value at a point \cite[Chap. 1]{AMS} and to the definition of  a distribution's value on a line  \cite{LVF0}.
\begin{definition}
	\label{defpoint}
	Let $s=1,2$. We say that a distribution $f\in H_\zco^{-s}$ has the
	value $f_0\in H^{-s}(\mathbb{R})$ on the line $x_1=0$, i.e. $f\big(0^+,(\cdot)_{[2]}\big)= f_0\big((\cdot)_{[2]}\big)$, if for each $\varphi\in H^s(\mathbb{R})$ and $\psi\in \mathscr D(\mathbb R_+)$ we have
	\begin{equation}
	\label{value}
	\left\langle\left\langle f(\alpha (\cdot)_{[1]},(\cdot)_{[2]}),\varphi((\cdot)_{[2]})\right\rangle_{[2]},\psi((\cdot)_{[1]})\right\rangle_{[1]}\rightarrow\left\langle\left\langle f_0,\varphi\right\rangle_{[2]},\psi\right\rangle_{[1]}\quad\text{as }\alpha\rightarrow 0^+, 
	\end{equation}
	where $\langle h(\alpha(\cdot)),\psi\rangle
	=\left\langle h((\cdot)),\frac1\alpha\psi\left(\frac{(\cdot)}{\alpha}\right)\right\rangle $ for $h\in H^{-s}(\mathbb R)$.
\end{definition}

\begin{remark}
	Let $\varphi \in H_\zco^s$, $s=\overline{0,2}$. Let $\widehat\varphi$ be 
	its even extension with respect to $x_1$, i.e.,
	$\widehat\varphi(x_1,x_2)=\varphi(x_1,x_2)$ if $x_1\ge0$ and  $\widehat\varphi(x_1,x_2)=\varphi(-x_1,x_2)$ if $x_1<0$, $x_2\in \mathbb R$. Then $\widehat\varphi \in \widehat H^s\left(\mathbb{R}^2\right)$, $s=\overline{0,2}$.
	The converse assertion is true for $s=0,1$, and it is not true for $s=2$. That is why the even extension with respect to $x_1$ of a distribution $f\in H_\zco^{-2}$
	may not belong to $\widehat{H}^{-2}\left(\mathbb{R}^2\right)$. However, the following theorem holds.
\end{remark}
\begin{theorem}
	\label{todext}
	Let $f\in H_\zco^0$ and there exist $f_{x_1}\big(0^+,(\cdot)_{[2]}\big)\in H^{-1}(\mathbb{R})$. Then $f_{x_1x_1}\in
	H_\zco^{-2}$ can be extended to a distribution $F\in \widehat{H}^{-2}\left(\mathbb{R}^2\right)$ such that $F$ is even with respect to $x_1$. This distribution is given by the formula
	\begin{equation}
	\label{oddext}
	F=\widehat f_{x_1x_1}-2f_{x_1}(0^+,\big(\cdot)_{[2]}\big)\delta_{[1]},
	\end{equation}
	where $\widehat f$ is the even extension of $f$ with respect to $x_1$.
\end{theorem}
In the case $f\in H_\zco^{1/2}$, corresponding theorem has been proved in \cite{LVF0}. The proof
of Theorem \ref{todext} is analogous to the proof of the mentioned theorem.

\section{Problem formulation and main results}\label{sect3}

We consider control system \eqref{eq}--\eqref{ic} in $H_\zco^{-l}$, $l=\overline{0,2}$, i.e. $\left(\frac
d{dt}\right)^s w:[0,T]\to H_\zco^{-2s}$, $s=0,1$, $w^0\in H_\zco^0$. We treat equality \eqref{bc} as the value of the distribution $w$ at $x_1=0$ with regard to Definition \ref{defpoint}.

Let $w^0,w(\cdot,t)\in H_\zco^0$, $t\in [0,T]$. Let $W^0$ and $W(\cdot,t)$ be the even extensions of $w^0$ and $w(\cdot,t)$ with
respect to $x_1$, respectively, $t\in [0,T]$. 
Consider the control system
\begin{align}
&W_t=\bigtriangleup W-2\delta u(t),&&  t\in(0,T),\label{eq1}
\\
&W\big((\cdot)_{[1]},(\cdot)_{[2]},0\big)=W^0, \label{ic1}
\end{align}
where $\left(\frac
d{dt}\right)^s W:[0,T]\to \widehat{H}^{-2s}\left(\mathbb{R}^2\right)$, $s=0,1$, $W^0\in \widehat{H}^0\left(\mathbb{R}^2\right)$, $\delta$ is the Dirac distribution in $\mathscr S'\left(\mathbb{R}^2\right)$. 

\begin{theorem}
	\label{lemequiv}
	\mbox{ }
	\begin{enumerate}
		\item \label{lemequiv-1}
		Let $w^0\in H_\zco^0$.
		If	$w$ is a solution to control system \eqref{eq}--\eqref{ic}, then $W$, its even extension with
		respect to $x_1$, is a solution to control system \eqref{eq1}, \eqref{ic1}.
		\item \label{lemequiv-2}
		Let $W^0\in \widehat{H}^0\left(\mathbb{R}^2\right)$.
		If	$W$ is a solution to control system \eqref{eq1}, \eqref{ic1},  then $w$, its restriction to $\mathbb R_+\times \mathbb R\times [0,T]$, is a solution to control system \eqref{eq}--\eqref{ic}.
	\end{enumerate}
\end{theorem}

\begin{proof}
	\ref{lemequiv-1}.
	Let $w$ be a solution to control system \eqref{eq}--\eqref{ic}. 
	According to Theorem \ref{todext}, $W$ is a solution to control system \eqref{eq1}, \eqref{ic1}.
	
	\ref{lemequiv-2}.
	Let $W$ be a solution to \eqref{eq1}, \eqref{ic1}. Let $w^0$ and $w(\cdot,t)$ be the restrictions of $W^0$ and $W(\cdot,t)$ to $\mathbb R_+\times \mathbb R$, respectively, $t\in [0,T]$.  According to Lemma \ref{lempoint} (see below),
	\begin{equation}
	\label{bc1}
	W_{x_1}\big(0^+,(\cdot)_{[2]},t\big)=\delta_{[2]}u(t)\quad\text{for almost all}\ t\in(0,T].
	\end{equation}
	Therefore, $w$ is a solution to \eqref{eq}--\eqref{ic}. 
\end{proof}

Due to Theorem \ref{lemequiv}, control systems \eqref{eq}--\eqref{ic} and \eqref{eq1}, \eqref{ic1} are equivalent. Therefore, basing on this reason, we will further consider control system  \eqref{eq1}, \eqref{ic1} instead of original system \eqref{eq}--\eqref{ic}.


Let $T>0$, $W^0\in\widehat{H}^0\left(\mathbb{R}^2\right)$.
By $\EuScript R_T\left(W^0\right)$, denote the set of all states $W^T\in\widehat{H}^0\left(\mathbb{R}^2\right)$ for which there exists a control $u\in L^\infty(0,T)$ such that there exists a unique solution $W$ to system \eqref{eq1}, \eqref{ic1} such that  $W\big((\cdot)_{[1]},(\cdot)_{[2]},T\big)=W^T$. 
\begin{definition}
	\label{def-cntr}
	A state $W^0\in\widehat{H}^0\left(\mathbb{R}^2\right)$ is said to be controllable to a target state $W^T\in\widehat{H}^0\left(\mathbb{R}^2\right)$ in a given time $T>0$ if
	$W^T\in\EuScript R_T\left(W^0\right)$.
\end{definition}
In other words, a state
$W^0\in\widehat{H}^0\left(\mathbb{R}^2\right)$ is said to be controllable to a target state $W^T\in\widehat{H}^0\left(\mathbb{R}^2\right)$ in a given time $T>0$ if there exists a control $u\in L^\infty(0,T)$ such that there exists a unique solution $W$ to system \eqref{eq1}, \eqref{ic1} and $W\big((\cdot)_{[1]},(\cdot)_{[2]},T\big)=W^T$.
\begin{definition}
	\label{def-a-cntr}
	A state $W^0\in\widehat{H}^0\left(\mathbb{R}^2\right)$ is said to be approximately controllable to a target state $W^T\in\widehat{H}^0\left(\mathbb{R}^2\right)$ in a given time $T>0$ if
	$W^T\in\overline{\EuScript R_T\left(W^0\right)}$, where the closure is considered in the space $\widehat{H}^0\left(\mathbb{R}^2\right)$.
\end{definition}
In other words, a state $W^0\in\widehat{H}^0\left(\mathbb{R}^2\right)$ is approximately
controllable  to a target state $W^T\in\widehat{H}^0\left(\mathbb{R}^2\right)$ in a given time
$T>0$ if for each $\varepsilon>0$, there exists a control $u_\varepsilon\in
L^\infty(0,T)$ such that there exists a unique solution $W_\varepsilon$ to
system \eqref{eq1}, \eqref{ic1} with $u=u_\varepsilon$ and $\left\|
W_\varepsilon\big((\cdot)_{[1]},(\cdot)_{[2]},T\big)-W^T \right\|^0<\varepsilon$.

Using the Poisson integral (see, e.g., \cite{VSV}), we obtain the unique solution to system \eqref{eq1},  \eqref{ic1} 
\begin{equation}
\label{sol11}
W(x,t)=\EuScript W_0(x,t)+\EuScript W_u(x,t), \quad x\in\mathbb{R}^2,\ t\in[0,T],
\end{equation}
where
\begin{align}
\label{sol1w0}
\EuScript W_0(x,t)&=\frac{1}{4\pi t}e^{-\frac{|x|^2}{4t}}*W^0(x),&&x\in\mathbb{R}^2,\ t\in[0,T],
\\
\label{sol1wu}
\EuScript W_u(x,t)&=-\frac{1}{\pi}\int_0^t \frac{1}{2\xi}e^{-\frac{|x|^2}{4\xi}} u(t-\xi)\,d\xi, && x\in\mathbb{R}^2,\ t\in[0,T].
\end{align}

Set $U_T^L=\{v\in L^\infty(0,T)\mid \|v\|_{L^\infty(0,T)}\leq L\}$ for $L>0$ and $T>0$.

According to \eqref{sol11}, we have
\begin{align}
\label{re0}
\kern-1.4ex\EuScript R_T\big(W^0\big)&\!=\!\left\{ W^T\in\widehat{H}^0\left(\mathbb{R}^2\right)\, \middle|\, \exists u\in L^\infty(0,T)\ \ 
W^T=\EuScript W_0(\cdot,T)+\EuScript W_u(\cdot,T)\right\},
\end{align}
in particular,
\begin{align}
\label{re1}
\EuScript R_T(0)&=\left\{ W^T\in\widehat{H}^0\left(\mathbb{R}^2\right)\, \middle|\, \exists u\in L^\infty(0,T)\ \ W^T=\EuScript W_u(\cdot,T)\right\}.
\end{align}
Denote also
\begin{align}
\label{re2a}
\EuScript R_T^L\big(W^0\big)&=\left\{ W^T\in\widehat{H}^0\left(\mathbb{R}^2\right)\, \middle|\, \exists u\in U_T^L\ \ 
W^T=\EuScript W_0(\cdot,T)+\EuScript W_u(\cdot,T)\right\},
\\
\label{re2}
\EuScript R_T^L(0)&=\left\{ W^T\in\widehat{H}^0\left(\mathbb{R}^2\right)\, \middle|\, \exists u\in U_T^L \ \
W^T=\EuScript W_u(\cdot,T)\right\}.
\end{align}

We obtain the following properties of a solution to system \eqref{eq1},  \eqref{ic1} 
\begin{theorem}
	\label{th-sol}
	Let $u\in L^\infty(0,T)$, $W^0\in\widehat{H}^0\left(\mathbb{R}^2\right)$. Then,
	\newcounter{enu}
	\renewcommand{\theenu}{(\roman{enu})}
	\renewcommand{\theenumi}{\upshape\thetheorem\ (\roman{enumi})}
	\begin{list}{\upshape{(\roman{enumi})}\hfill}{\usecounter{enumi} \leftmargin=5ex \labelwidth=5ex \labelsep=0ex \topsep=1ex \listparindent=1.5ex}	
		\item \label{th-sol-i} \refstepcounter{enu} \label{th-sol-ia}
		$\EuScript W_0(\cdot,t)\in\widehat H^0(\mathbb R^2)$, $t\in[0,T];$
		\item \label{th-sol-ii}  \refstepcounter{enu} \label{th-sol-iia}
		$\EuScript W_0(\cdot,t)\in C^\infty(\mathbb R^2)$, $t\in(0,T];$
		\item \label{th-sol-iii}  \refstepcounter{enu} \label{th-sol-iiia}
		if $W^0\in\EuScript H$, then $\EuScript W_0(\cdot,t)\in\EuScript H$, $t\in[0,T];$
		\item \label{th-sol-iv}  \refstepcounter{enu} \label{th-sol-iva}
		$\EuScript W_u(\cdot,t)\in\EuScript H$ and $\|\EuScript W_u(\cdot,t)\|^0\leq \frac{2}{\sqrt{\pi}}(t+1) \|u\|_{L^\infty(0,T)}$, $t\in[0,T].$
	\end{list}
\end{theorem}
The \hyperlink{th34}{proof of the theorem} is given in Section \ref{sect6}.

With regard to Theorem \ref{th-sol}, we get the following properties of the sets $\EuScript R_T(g)$ and $\EuScript R_T^L(g)$.
\begin{theorem}
	\label{reachprop}
	Let $T>0$, $g\in \widehat{H}^0\left(\mathbb{R}^2\right)$. We have
	\renewcommand{\theenumi}{\upshape\thetheorem\ (\roman{enumi})}
	\begin{list}{\upshape{(\roman{enumi})}\hfill}{\usecounter{enumi} \leftmargin=5ex \labelwidth=5ex \labelsep=0ex \topsep=1ex \listparindent=1.5ex}
		\item\label{laa1}
		$\displaystyle\EuScript R_T(0)=\bigcup_{L>0}\EuScript R_T^L(0)\subset \EuScript H;$
		\item\label{laa2}
		$\displaystyle\EuScript R_T^L(0)\subset\EuScript R_T^{L'}(0)$, $0<L<L';$
		\item\label{laa3}
		$\displaystyle f\in\EuScript R_T^1(0)\Leftrightarrow Lf\in\EuScript R_T^L(0)$, $L>0;$
		\item\label{laa4}
		$\displaystyle f\in\EuScript R_T^L(g) \Leftrightarrow \left(f-\frac{1}{4\pi T}e^{-\frac{|\cdot|^2}{4T}}*g\right)\in\EuScript R_T^L(0)$, $L>0;$
		\item\label{laa5}
		$\displaystyle f\in\EuScript R_T(g) \Leftrightarrow \left(f-\frac{1}{4\pi T}e^{-\frac{|\cdot|^2}{4T}}*g\right)\in\EuScript R_T(0)$.
	\end{list}
\end{theorem}

Consider also the sets
\begin{equation}
\label{r}
\RR=\Psi \EuScript R_T(0)\quad \text{and}\quad \RRL=\Psi \EuScript R_T^L(0),
\end{equation}
where Theorem \ref{th-sol-iv} is taken into account. 
Put $\EuScript Y_u(\cdot,t)=\Psi \EuScript W_u(\cdot,t)$, $t\in[0,T]$.  Then, we have
\begin{equation}
\label{u}
\EuScript Y_u(r,t)=-\frac1\pi \int_0^t \frac1{2\xi} e^{-\frac r{4\xi}} u(t-\xi)\, d\xi, \quad r\in\mathbb R_+.
\end{equation}
With regard to \eqref{re1} and \eqref{re2}, we obtain
\begin{align}
\label{r1}
\RR&=\left\{ Y\in \LL \ \middle|\ \exists u\in L^\infty(0,T)\quad Y=\EuScript Y_u(\cdot,T)\right\},
\\
\label{r2}
\RRL&=\left\{Y\in\LL\ \middle|\ \exists u\in U_T^L \quad 
Y=\EuScript Y_u(\cdot,T)\right\}.
\end{align}
Since $\Psi$ is an isomorphism of $\HH$ and $\LL$ (see Proposition \ref{prop-1-i}), we have
\begin{equation}
\label{rrr}
\overline{\RR}=\Psi \overline{\EuScript R_T(0)}\quad \text{and}\quad \overline{\RRL}=\Psi \overline{\EuScript R_T^L(0)}.
\end{equation}

Properties of the sets $\RR$, $\overline{\RR}$, $\RRL$, and $\overline{\RRL}$ are studied in Section \ref{prop}. By using these properties, we obtain the main results of the paper.


\subsection{Controllability under controls bounded by a hard constant.}

\label{sect4}

Let us find conditions under which an initial state $W^0\in\widehat{H}^0\left(\mathbb{R}^2\right)$ is controllable to  a target state $W^T\in\widehat{H}^0\left(\mathbb{R}^2\right)$ in a given time $T>0$. 

First, consider  necessary conditions for $f\in\EuScript R_T^L(0)$.
\begin{theorem} 
	\label{thnec} 
	Let $L>0$ and $T>0$. If $f\in\EuScript R_T^L(0)$, then  $f\in\HH$ and we have 
	\begin{equation}
	\label{nec}
	|f(x)|\leq \frac L{2\pi} e^{-\frac{|x|^2}{4T}}\ln \left(1+\frac{4T}{|x|^2}\right),\quad x\in\mathbb{R}^2\setminus\{0\}.
	\end{equation}
\end{theorem}

The \hyperlink{th36}{proof of the theorem} is given in Section \ref{sect6}.

\begin{theorem} 
	\label{entire}
	Let $L>0$ and $T>0$. Let also $f\in\EuScript R_T^L(0)$, $F=\mathscr F f$ and $G=\Psi F$. Then $G$ can be extended to an entire function $G_e$  of   order $\leq 1$ and   type $\leq T$. Moreover, $F$ can be also extended to an entire function $F_e$ and
	\begin{equation}
	\label{fff3}
	F_e(s)=G_e\big(s_1^2+s_2^2\big),\quad s=(s_1,s_1)\in\mathbb C^2,
	\end{equation}
	and $F_e$ is of   order $\leq 2$ and   type $\leq T$. In addition,
	\begin{align}
	\label{ent00}
	\big|F_e(s)\big|\leq \big|G_e\big(s_1^2+s_2^2\big)\big|\leq \frac L\pi \frac{e^{T|s|^2}-1}{|s|^2},\quad s\in \mathbb C^2.
	\end{align}
\end{theorem}

The \hyperlink{th37}{proof of the theorem} is given in Section \ref{sect6}.

According to Example \ref{ex-nec} below, condition \eqref{nec} is only necessary for $f\in\EuScript R_T^L(0)$, but it is not sufficient. However, if $f$ satisfies \eqref{nec}, its Fourier transform can be extended to an entire function of   order $\leq2$ and   type $\leq T$ (cf. Theorem \ref{entire}).

\begin{theorem} 
	\label{cond-f}
	Let $L>0$, $T>0$, $f\in \EuScript H$, and condition \eqref{nec} hold for $f$. Let also  $F=\mathscr F f$ and $G=\Psi F$. Then $G$ can be extended to an entire function $G_e$  of   order $\leq 1$ and   type $\leq T$. Moreover, $F$ can be also extended to an entire function $F_e$, the extension $F_e$ is given by \eqref{fff3}, and $F_e$ is of order $\leq 2$ and type $\leq T$.
\end{theorem}

The \hyperlink{th38}{proof of the theorem} is given in Section \ref{sect6}.

Thus, condition \eqref{nec} is not sufficient for $f\in\EuScript R_T^L(0)$, but it  guarantees the necessary condition from Theorem \ref{entire} holds for $\mathscr F f$.

Now, we consider a necessary and sufficient condition  for controllability in a given time $T>0$ under controls bounded by a hard constant.
Denote
\begin{equation}
\label{www} 
W^T_0=W^T-\EuScript W_0(\cdot,T).
\end{equation}

\begin{theorem}
	\label{thmom} 
	Let $L>0$, $T>0$, $W^0\in\widehat{H}^0\left(\mathbb{R}^2\right)$, $W^T\in\widehat{H}^0\left(\mathbb{R}^2\right)$. Let also $W^T_0\in \EuScript H$, condition \eqref{nec} hold for $W^T_0$, and
	\begin{equation} 
	\label{moment}
	\omega_n=-2\frac{n!}{(2n)!}\int_0^\infty\int_0^\infty x_1^{2n}W^T_0(x_1,x_2)\,dx_1\,dx_2,\qquad n\in \mathbb N_0.
	\end{equation}
	Then $W^T\in\EuScript R_T^L\left(W^0\right)$ iff 
	\begin{equation}
	\label{moment1}
	\exists u \in U_T^L\ \forall n\in\mathbb N_0\quad 
	\int_0^T \xi^n u(T-\xi)\,d\xi=\omega_n.
	\end{equation}
\end{theorem}

The \hyperlink{th39}{proof of the theorem} is given in Section \ref{sect6}.

Taking into account Definition \ref{def-cntr}, we get
\begin{corollary}
	\label{thmom-cor} 
	Let $L>0$, $T>0$, $W^0\in\widehat{H}^0\left(\mathbb{R}^2\right)$, $W^T\in\widehat{H}^0\left(\mathbb{R}^2\right)$. Let also $W^T_0\in \EuScript H$, condition \eqref{nec} hold  for $W^T_0$, and $\{\omega_n\}_{n=1}^\infty$ be determined by \eqref{moment}.
	Then the state $W^0\in\widehat{H}^0\left(\mathbb{R}^2\right)$ is  controllable to the target state $W^T\in\widehat{H}^0\left(\mathbb{R}^2\right)$ in a given time $T>0$ iff \eqref{moment1} holds.
\end{corollary}

Now, we consider a  sufficient condition  for approximate controllability in a given time $T>0$ under controls bounded by a hard constant.

\begin{theorem} 
	\label{thmomap}  Let $L>0$, $T>0$, $W^0\in\widehat{H}^0\left(\mathbb{R}^2\right)$, $W^T\in\widehat{H}^0\left(\mathbb{R}^2\right)$, $W^T_0\in\EuScript H$, and condition \eqref{nec} hold for
	$W^T_0$. Let
	$\{\omega_n\}_{n=0}^\infty$ be defined by \eqref{moment}. If for each
	$N\in\mathbb{N}$ there exists $u_N\in U_T^L$ such that 
	\begin{equation}
	\label{moment9}
	\int_0^T\xi^nu_N(T-\xi)d\xi=\omega_n,\quad n=\overline{0,N},
	\end{equation}
	then $W^T\in\overline{\EuScript R_T^L\left(W^0\right)}$, where the closure is considered in $\widehat{H}^0\left(\mathbb{R}^2\right)$.
\end{theorem}

The \hyperlink{th311}{proof of the theorem} is given in Section \ref{sect6}.

Taking into account Definition \ref{def-a-cntr}, we get
\begin{corollary}
	\label{thmomap-cor}  
	Let $L>0$, $T>0$, $W^0\in\widehat{H}^0\left(\mathbb{R}^2\right)$, $W^T\in\widehat{H}^0\left(\mathbb{R}^2\right)$, $W^T_0\in\EuScript H$, and condition \eqref{nec} hold for
	$W^T_0$. Let
	$\{\omega_n\}_{n=0}^\infty$ be defined by \eqref{moment}. If for each
	$N\in\mathbb{N}$ there exists $u_N\in U_T^L$ such that 
	\eqref{moment9} holds,
	then the state $W^0\in\widehat{H}^0\left(\mathbb{R}^2\right)$ is approximately controllable to the target state $W^T\in\widehat{H}^0\left(\mathbb{R}^2\right)$ in a given time $T>0$.
\end{corollary}

We can see that the controllability problems were reduced to the Markov power moment problems in Theorems \ref{thmom} and \ref{thmomap}. These Markov power moment problems may be solved by using the algorithms given in \cite{LGOR, KS3}. Similar results were obtained for controllability problems for the heat equation on a half-axis \cite{FKh, FKh2} and on a half-plane \cite{FKh3}. However, the description of the set $\EuScript R_T^L\left(W^0\right)$ is given in principally different way in the present paper (see Theorems \ref{thnec}, \ref{entire}, \ref{cond-f}).  As a result, the necessary condition \eqref{nec} essentially differs from the necessary conditions obtained in the mentioned papers:
it is given in the form of an estimate for a function belonging to $\EuScript R_T^L\left(W^0\right)$ in contrast to the conditions in the form of estimates for integrals with special weights of a such function in 
\cite{FKh, FKh2, FKh3}. In addition, as a consequence of the different necessary condition \eqref{nec}, the proofs of Theorems \ref{thmom} and \ref{thmomap} also differ from their analogues in the mentioned papers.



\subsection{Approximate controllability.}
\label{sect5}

Consider the problem of approximate controllability for system \eqref{eq1}, \eqref{ic1} under controls from $L^\infty(0,T)$ unlike Subsection \ref{sect4}, where we consider this system under controls bounded by a hard constant. We have the following main theorem.

\begin{theorem}
	\label{thclos}
	Let $T>0$. We have $\overline{\EuScript R_T(0)}=\EuScript{H}$.
\end{theorem} 
The \hyperlink{th313}{proof of the theorem} is given in Section \ref{sect6}.
This theorem yields
\begin{theorem}
	\label{thappc}
	Let $T>0$. A state $W^0\in\widehat{H}^0\left(\mathbb{R}^2\right)$ is approximately controllable to a state $W^T\in\widehat{H}^0\left(\mathbb{R}^2\right)$
	in a given time $T$ iff $W^T_0\in\EuScript H$.
\end{theorem}


\subsection{Lack of controllability to the origin.}
\label{sect5a}

For $W^0\in\widehat{H}^0\left(\mathbb{R}^2\right)$ and $W^T\in\widehat{H}^0\left(\mathbb{R}^2\right)$ we have $W^T\in\overline{\EuScript R_T\big(W^0\big)}$ iff  $W^T_0\in\EuScript H$ according to Theorem \ref{thappc}.
However, $0\notin \EuScript R_T\big(W^0\big)$ for all nonzero $W^0\in\widehat{H}^0\left(\mathbb{R}^2\right)$, i.e. the following theorem holds.
\begin{theorem}
	\label{thnulcontr}
	If a state $W^0\in\widehat{H}^0\left(\mathbb{R}^2\right)$ is controllable to the state $W^T=0$ in a given time $T>0$, then $W^0=0$.
\end{theorem}

The \hyperlink{th315}{proof of the theorem} is given in Section \ref{sect6}.

\section{Properties of the sets \texorpdfstring{$\RR$}{R T} and \texorpdfstring{$\RRL$}{R T L}}
\label{prop}

First, consider  necessary conditions for $g\in\RRL$.
\begin{theorem} 
	\label{pthnec} 
	Let $T>0$, $L>0$ and $g\in\RRL$. Then
	\begin{equation}
	\label{pnec}
	\frac{ e^{\frac{(\cdot)}{4T}} g}{\ln \left(1+\frac{4T}{(\cdot)}\right)} \in L^\infty(\mathbb R_+) .
	\end{equation}
	In addition, we have
	$$
	\left\|\frac{ e^{\frac{(\cdot)}{4T}} g}{\ln \left(1+\frac{4T}{(\cdot)}\right)} \right\|_{L^\infty(\mathbb R_+)}
	\leq \frac L{2\pi}.
	$$
\end{theorem}

\begin{proof}
	According to \eqref{r2}, there exists $u\in U_T^L$ such that
	$$
	g(r)=-\frac1\pi \int_0^T\frac{e^{-\frac r{4\xi}}}{2\xi} u(T-\xi)\, d\xi, \quad r\in \mathbb R_+.
	$$
	Setting $y=r/(4\xi)$ and taking into account \cite[5.1.1, 5.1.20]{HB}, we get
	\begin{equation*}
	|g(r)|\leq \frac L{2\pi} \int_{r/(4T)}^\infty \frac{e^{-y}}{y} \, dy= \frac L{2\pi} E_1\left(\frac{r}{4T}\right) \leq \frac L{2\pi} e^{-\frac r{4T}} \ln\left(1+\frac{4T}{r}\right),\quad r>0, 
	\end{equation*}
	where $E_1(\xi)=\int_{\xi}^\infty \big(e^{-t}/t\big)\, dt$, $\xi\in\mathbb R$. Therefore,  \eqref{pnec} holds and the estimate for the norm is true. 
\end{proof}

We need the following formula
\begin{equation}
\label{phi-e}
\Phi\left(e^{-\alpha(\cdot)}\right)=\frac1{2\alpha} e^{-\frac{(\cdot)}{4\alpha}},\quad \alpha\in \mathbb R_+.
\end{equation}
To prove it, set $\alpha\in \mathbb R_+$ and $q(r)=e^{-\alpha r}$, $r\in\mathbb R_+$. Expanding the Bessel function into the power series, we get
$$
\big(\Phi q\big)(\rho)
=\frac12\int_0^\infty e^{-\alpha r} J_0\big(\sqrt{r\rho}\big)\, dr
=\frac12\sum_{m=0}^\infty \frac{(-1)^m \rho^m}{(m!)^2 2^{2m}} \int_0^\infty r^m e^{-\alpha r}\, dr,\quad \rho\in \mathbb R_+.
$$
Since 
\begin{align}
\label{int-e}
\int_0^\infty r^m e^{-\alpha r}\, dr
&=(-1)^m \left(\frac d{d\alpha}\right)^m\int_0^\infty e^{-\alpha r}\, dr
\notag\\
&= (-1)^m \left(\frac d{d\alpha}\right)^m \frac1\alpha
=\frac{m!}{\alpha^{m+1}},
\end{align}
we get
$$
\big(\Phi q\big)(\rho)
=\frac1{2\alpha} \sum_{m=0}^\infty \frac{(-1)^m}{m!} \left(\frac{\rho}{4\alpha}\right)^m
=\frac1{2\alpha} e^{-\frac{\rho}{4\alpha}},\quad \rho\in \mathbb R_+,
$$
i.e. \eqref{phi-e} holds.

\begin{theorem} 
	\label{pentire}
	Let $T>0$, $L>0$, $g\in\RRL$, and $G=\Phi g$.   Then $G$ can be extended to an entire function $G_e$  of   order $\leq 1$ and   type $\leq T$ and
	\begin{equation}
	\label{ent0}
	|G_e(z)|\leq \frac L\pi \frac{e^{T|z|}-1}{|z|},\quad z\in \mathbb C.
	\end{equation}
\end{theorem}

\begin{proof}
	Since $g\in\RRL$, there is $u\in U_T^L$ such that $g=\EuScript Y_u(\cdot,T)$ according to \eqref{r2}. With regard to \eqref{phi-e}, we get
	$$
	G(\rho)=-\frac1\pi \int_0^T e^{-\xi \rho} u(T-\xi)\, d\xi, \quad \rho\in \mathbb R_+.
	$$
	Hence $G$ can be extended to the entire function $G_e$ by the formula
	\begin{equation}
	\label{ent2}
	G_e(z)=-\frac1\pi \int_0^T e^{-\xi z} u(T-\xi)\, d\xi, \quad z\in \mathbb C.
	\end{equation}
	Evidently,
	\begin{equation*}
	|G_e(z)|\leq \frac L\pi \int_0^T e^{\xi |z|}\, d\xi
	= \frac L\pi  \frac{e^{T|z|}-1}{|z|},\quad z\in \mathbb C.
	\end{equation*}
	The theorem is proved.
\end{proof}


\begin{example}
	\label{ex-nec}
	Let  $T>0$ and
	$$
	g(r)=- \frac2{\pi T} e^{-\frac{r}{2T}},\quad r\geq 0.
	$$	
	Let us verify condition \eqref{pnec}  for $g$. Put $v(\xi)=2\ln(1+1/\xi)-e^{-\xi}$, $\xi>0$. Then $v'(\xi)= -2/(\xi+\xi^2)+e^{-\xi}$, $\xi>0$. Since $\xi+\xi^2 < 2+2\xi+\xi^2<2e^\xi$, $\xi>0$, we have $v'(\xi)<0$, $\xi>0$, i.e. $v$ decreases on $(0,+\infty)$. We have $v(\xi)\to+\infty$ as $\xi\to0^+$ and $v(\xi)\to0$ as $\xi\to+\infty$. Therefore, $v(\xi)>0$, $\xi>0$, i.e.
	$$
	e^{-\xi}\leq 2 \ln\left(1+\frac1\xi\right),\quad \xi>0.
	$$
	Setting $\xi=r/(4T)$ and applying this estimate to $g$, we get
	$$
	|g(r)|\leq \frac{4}{\pi T} e^{-\frac{r}{4T}} \ln\left(1+
	\frac{4T}r\right),\quad r>0.
	$$
	Therefore, condition  \eqref{pnec} holds for $g$.
	
	Let us try to find $u\in U_T^L$ such that
	$$
	g(r)=\EuScript Y_u(x,T)=-\frac1\pi \int_0^T \frac1{2\xi} e^{-\frac{r}{4\xi}} u(T-\xi)\, d\xi,\quad r\in \mathbb R_+.
	$$
	Applying the operator $\Phi$, we get
	\begin{align}
	\label{est-ex}
	-\frac2\pi e^{-T\rho/2}&=\big(\Phi g \big)(\rho) =-\frac1\pi\int_0^T e^{-\xi\rho} u(T-\xi)\, d\xi
	\notag\\
	&=-\sqrt{\frac2\pi} \big(\mathscr F \mathcal U_T\big)(-i\rho), \quad \rho\in\mathbb R_+,
	\end{align}
	where
	$$
	\mathcal U_T(\xi)=\begin{cases}
	u(T-\xi),& \xi \in[0,T],\\
	0, & \xi\in\mathbb R \setminus [0,T].
	\end{cases}
	$$
	Due to the Paley--Wiener theorem, $\mathscr F \mathcal U_T$ can be extended to an entire function. Replacing $\rho$ by $i\mu$, we obtain
	\begin{equation}
	\label{est-ex-1}
	\sqrt{\frac2\pi}e^{-iT\mu/2}=\big(\mathscr F \mathcal U_T\big)(\mu), \quad \mu \in \mathbb C.
	\end{equation}
	Therefore, $\mathcal U_T(\xi)=2\delta(\xi-T/2)$ is the unique solution to equation \eqref{est-ex-1}. Hence $u(\xi)=2\delta(\xi-T/2)$ is the unique solution to equation \eqref{est-ex}. But this function $u$ does not belong to $L^\infty(0,T)$. Therefore, $g\notin \RRL$ for any $T>0$ and $L>0$ although condition \eqref{pnec} holds for it.
\end{example}

Thus, condition \eqref{pnec} is only necessary for $g\in \RRL$, but it is not sufficient. However, if $g$ satisfies \eqref{pnec}, $\Phi g$ can be extended to an entire function of   order $\leq1$ and   type $\leq T$ (cf. Theorem \ref{pentire}).

\begin{theorem} 
	\label{pcond-f}
	Let $T>0$,  $g\in\LL$,   $G=\Phi g$, and condition \eqref{pnec} hold for $g$. Then $G$ can be extended to an entire function $G_e$  of   order $\leq 1$ and   type $\leq T$. 
\end{theorem}

\begin{proof}
	According to Proposition \ref{prop-2-ii}, we have
	\begin{equation}
	\label{fff00}
	G(\rho)=(\Phi g)(\rho)
	=\frac12 \lim_{N\to\infty} \int_0^N g(r) J_0\big(\sqrt{r\rho}\big)\, dr,\quad \rho \in \mathbb R_+.
	\end{equation}
	Setting $M=\left\|e^{\frac{(\cdot)}{4T}} g\middle/ \ln \left(1+\frac{4T}{(\cdot)}\right) \right\|_{L^\infty(\mathbb R_+)}$, we obtain from \eqref{pnec} that
	\begin{equation}
	\label{fff0}
	|g(r)|\leq M e^{-\frac{r}{4T}} \ln \left(1+\frac{4T}{r}\right)
	\leq 2\sqrt{2T} M \frac{ e^{-\frac{r}{4T}}}{\sqrt r}, \quad r \in \mathbb R_+,
	\end{equation}
	where the estimate 
	\begin{equation}
	\label{ln}
	\ln\big(1+y^2\big)\leq \sqrt2 y,\quad y\in \mathbb R_+,
	\end{equation}
	obtained from the obvious estimate
	$$
	1+y^2\leq 1+\sqrt 2 y + y^2 \leq e^{\sqrt2 y},\quad y\in \mathbb R_+,
	$$
	has been also used. 
	It follows from   \eqref{fff0} that 
	\begin{equation}
	\label{fff11}
	G_e(z)=\frac12 \int_0^\infty g(r) J_0\big(\sqrt{r z}\big)\, dr,\quad z \in \mathbb C,
	\end{equation}
	is an entire function because $J_0$ is an entire function. Due to \eqref{fff00}, we have
	\begin{equation}
	\label{fff22}
	G_e(r)=G(r),\quad r \in \mathbb R_+.
	\end{equation}
	Taking into account \eqref{fff0}, we get
	\begin{align}
	\label{fff2}
	|G_e(z)|
	&\leq \sqrt{2T} M 
	\int_0^\infty  \big|J_0\big(\sqrt{r z}\big)\big| \frac{e^{-\frac{r}{4T}}}{\sqrt r}\, dr
	=2 \sqrt{2T} M 
	\int_0^\infty  \big|J_0\big(y\sqrt{z}\big)\big| e^{-\frac{y^2}{4T}}\, dy
	\notag\\
	&
	\leq  2\sqrt{2T} M \sum_{m=0}^\infty \frac{|z|^m}{(m!)^2 2^{2m}} \int_0^\infty y^{2m} e^{-\frac{y^2}{4T}}\, dy,\quad z \in \mathbb C.
	\end{align}
	Let us calculate the last integral. We have
	\begin{align}
	\label{exp}
	\int_0^\infty y^{2m} e^{-\alpha y^2}\, dy
	&=(-1)^m \left(\frac{d}{d\alpha}\right)^m \int_0^\infty  e^{-\alpha y^2}\, dy
	= \frac{(-1)^m}2 \left(\frac{d}{d\alpha}\right)^m \sqrt{\frac{\pi}{\alpha}}
	\notag\\
	&
	=\frac{\sqrt\pi}{2} \frac{(2m-1)!!}{2^m \alpha^{m+1/2}}, \quad \alpha\in \mathbb R_+,\ m\in\mathbb N_0,
	\end{align}
	where we set $(-1)!!=1$.
	Setting $\alpha=1/(4T)$ and continuing \eqref{fff2}, we obtain
	\begin{align*}
	|G_e(z)|
	&\leq  2\sqrt{2\pi}T M \sum_{m=0}^\infty \frac{(4T)^m|z|^m}{m!2^{2m}} \frac{(2m-1)!!}{(2m)!!}
	\\
	&\leq 2\sqrt{2\pi}T M \sum_{m=0}^\infty \frac{T^m|z|^m}{m!} 
	\leq 4\sqrt{2}T M e^{T|z|},\quad z \in \mathbb C. 
	\end{align*}
	Thus, $G_e$ is an entire function of order $\leq1$ and type $\leq T$, and condition \eqref{fff22} holds for it.
\end{proof}

Thus, condition \eqref{pnec} is not sufficient for $g\in\RRL$, but it  guarantees the necessary condition from Theorem \ref{pentire} holds for $\Phi g$.


\begin{theorem} 
	\label{pthmom} 
	Let $L>0$, $T>0$, $g\in\LL$, and condition \eqref{pnec} hold for $g$. Let also 
	\begin{equation}
	\label{pmoment}
	\gamma_n=-\frac\pi{2^{2n+1}n!}\int_0^\infty r^n g(r)\,dr,\quad n\in \mathbb N_0.
	\end{equation}
	Then  $g\in\RRL$ iff 
	\begin{equation}
	\label{pmoment1}
	\exists u \in U_T^L\ \forall n\in\mathbb N_0\quad 
	\int_0^T \xi^n u(T-\xi)\,d\xi=\gamma_n.
	\end{equation}
\end{theorem}

\begin{proof}
	Due to  \eqref{r2}, we have
	\begin{equation}
	\label{moment3}
	g\in \RRL
	\Leftrightarrow
	\big( \exists u \in U_T^L\quad  g=\EuScript Y_u(\cdot,T) \big).
	\end{equation}
	Put
	$G=\Phi g$, $\EuScript G_u(\cdot,T)=\Phi \EuScript Y_u(\cdot,T)$. According to Theorem \ref{pentire}, $G$ can be extended to an entire function $G_e$. Taking into account \eqref{fff11}, we  obtain
	\begin{align}
	\label{moment5}
	G_e(z)&=\frac12 \int_0^\infty g(r) J_0\big(\sqrt{r z}\big)\, dr
	\notag\\
	&=\frac12 \sum_{n=0}^\infty \frac{(-1)^n z^n}{(n!)^2 2^{2n}} \int_0^\infty r^n g(r)\, dr
	=-\frac1\pi \sum_{n=0}^\infty \frac{(-1)^n \gamma_n}{n!} z^n
	,\quad z \in \mathbb C.
	\end{align}
	Taking into account \eqref{u} and \eqref{phi-e}, we get
	\begin{align}
	\label{moment6}
	\big(\Phi \EuScript Y_u(\cdot,T)\big)(z)& 
	=-\frac1\pi \int_0^T e^{-\xi z} u(T-\xi)\, d \xi
	\notag\\
	&= -\frac{1}{\pi} \sum_{n=0}^\infty \frac{(-1)^n z^n}{n!} \int_0^T \xi^n u(T-\xi)\, d \xi ,\quad z\in\mathbb C.
	\end{align}
	It follows from \eqref{moment3}--\eqref{moment6} that $g\in\RRL$ iff  \eqref{pmoment1} holds.
\end{proof}


\begin{theorem} 
	\label{pthmomap}  
	Let $L>0$, $T>0$, $g\in\LL$, condition \eqref{pnec} hold for $g$, and
	$\{\gamma_n\}_{n=0}^\infty$ be defined by \eqref{pmoment}. If for each
	$N\in\mathbb{N}$ there exists $u_N\in U_T^L$ such that 
	\begin{equation}
	\label{pmoment9}
	\int_0^T\xi^nu_N(T-\xi)d\xi=\gamma_n,\quad n=\overline{0,N},
	\end{equation}
	then $g\in\overline{\RRL}$, where the closure is considered in $\LL$.
\end{theorem}

\begin{proof}
	Let $N\in \mathbb N$.
	Set $g_N=\EuScript Y_{u_N}(\cdot,T)$, $G=\Phi g$, and $G_N=\Phi g_N$.
	According to Theorem \ref{pcond-f}, we have
	$G$ can be extended to an entire function $G_e$  of   order $\leq1$ and   type $\leq T$. In addition, \eqref{moment5} holds for it.
	Setting 
	\begin{equation}
	\label{moment100}
	\gamma_n^N=\int_0^T\xi^nu_N(T-\xi)d\xi,\quad n\in \mathbb N_0,
	\end{equation}
	and taking into account  \eqref{moment6}, we get  
	\begin{align}
	\label{moment10}
	G_N(z)& 
	=-\frac1\pi \int_0^T e^{-\xi z} u_N(T-\xi)\, d \xi
	\notag\\
	&= -\frac{1}{\pi} \sum_{n=0}^\infty \frac{(-1)^n \gamma_n^N}{n!} z^n ,\quad z\in\mathbb C.
	\end{align}
	It follows from  \eqref{moment5}, \eqref{pmoment9},  and \eqref{moment10} that 
	\begin{equation}
	\label{moment13}
	G_e(z)-G_N(z)=-\frac1\pi \sum_{n=N+1}^\infty \frac{(-1)^n}{n!} \big(\gamma_n-\gamma_n^N\big) z^n,\quad z\in \mathbb C.
	\end{equation}
	
	Let $\varepsilon>0$ be fixed. Then there exists $A_\varepsilon>0$ such that
	\begin{equation}
	\label{moment15}
	\int_{A_\varepsilon}^\infty \left|G(\rho)-G_N(\rho) \right|^2 d\rho < \varepsilon^2.
	\end{equation}
	With regard to \eqref{moment13}, we get
	\begin{equation}
	\label{moment16}
	\left|G(\rho)-G_N(\rho) \right|
	\leq \frac1\pi \sum_{n=N+1}^\infty \frac{1}{n!} \big|\gamma_n-\gamma_n^N\big| A_\varepsilon^{n},\quad 
	\rho\in (0,A_\varepsilon].
	\end{equation}
	Now, let us estimate $\gamma_n$ and $\gamma_n^N$. It follows from \eqref{pnec} that 
	$$
	|g(r)|\leq M e^{-\frac{r}{4T}} \ln \left(1+\frac{4T}{r}\right), \quad r\in\mathbb R_+,
	$$
	for some $M>0$. Taking into account \eqref{ln} and \eqref{pmoment}, we get
	\begin{align*}
	|\gamma_n|&\leq \frac{M\pi}{2^{2n+1}n!}
	\int_0^\infty r^n e^{-\frac{r}{4T}} \ln \left(1+\frac{4T}{r}\right)\,dr
	\\
	&\leq \frac{\sqrt{2T}M\pi}{2^{2n-1}n!} \int_0^\infty r^n e^{-\frac{r}{4T}}\frac{dr}{2\sqrt r}
	=\frac{\sqrt{2T}M\pi}{2^{2n-1}n!} \int_0^\infty y^{2n} e^{-\frac{y^2}{4T}}\,dr,\quad n\in\mathbb N_0. 
	\end{align*}
	Taking into account \eqref{exp}, we obtain
	\begin{equation}
	\label{omega}
	|\gamma_n|\leq \frac{M\pi^{3/2}}{2^{2n+1/2}} \frac{(2n-1)!!}{(2n)!!}(4T)^{n+1} \leq M(2\pi)^{3/2} T^{n+1},\quad n\in\mathbb N_0.
	\end{equation}
	It follows from \eqref{moment100} that 
	\begin{equation}
	\label{omegaN}
	\big|\gamma_n^N\big|\leq L\int_0^T\xi^n\, d\xi=L\frac{T^{n+1}}{n+1},\quad n\in\mathbb N_0.
	\end{equation}
	According to \eqref{omega} and \eqref{omegaN}, we get
	$$
	\big|\gamma_n-\gamma_n^N\big|\leq \left(M(2\pi)^{3/2} +\frac{L}{n+1}\right)T^{n+1} \leq \pi C T^n,\quad n\in\mathbb N_0,
	$$
	where $C=T(M(2\pi)^{3/2} +L/(n+1))$.
	With regard to \eqref{moment16}, we have
	\begin{align*}
	\left|G(\rho)-G_N(\rho) \right|
	\leq C \sum_{n=N+1}^\infty \frac{\big(T A_\varepsilon \big)^n}{n!}
	&=C\left(e^{T A_\varepsilon } -\sum_{n=0}^N \frac{\big(T A_\varepsilon \big)^n}{n!}\right) 
	\\
	&\leq C e^{T A_\varepsilon } \frac{(T A_\varepsilon )^{N+1}}{(N+1)!},\quad 
	\rho\in (0,A_\varepsilon],
	\end{align*}
	because
	$$
	\left|e^y -\sum_{n=0}^N \frac{y^n}{n!}\right|
	\leq e^{|y|} \frac{|y|^{N+1}}{(N+1)!},\quad N\in\mathbb N,
	$$
	according to the Taylor formula. Therefore,
	$$
	\left(\int_0^{A_\varepsilon} \left|G(\rho)-G_N(\rho) \right|^2 d\rho\right)^{1/2}
	\leq  C \sqrt{A_\varepsilon} e^{T A_\varepsilon } \frac{(T A_\varepsilon )^{N+1}}{(N+1)!}\underset{N\to\infty}\to0.
	$$
	Taking into account \eqref{moment15}, we conclude that
	$$
	\left\|g- g_N\right\|_{\LL}
	=\left\| G- G_N\right\|_{\LL}\underset{N\to\infty}\to0.
	$$
	Therefore, $g\in\overline{\RRL}$.
\end{proof}


Put
\begin{align}
\label{hhh}
\varphi_{n}(\rho)&=\rho^{n}e^{-T\rho}, &\rho&\in\mathbb R_+, &n&\in\mathbb N_0, &&
\\
\label{hhhl}
\varphi_{n}^l(\rho)&=\rho^{n}e^{-T\rho} \left(\frac{e^{\rho/l}-1}{\rho/l}\right)^{n+1}, & \rho&\in\mathbb R_+, &n&\in\mathbb N_0,& l&\in\mathbb N.
\end{align}

First we consider the system  $\{\varphi_n\}_{n=0}^\infty$. It is well-known that it is complete in $\LL$.

The following lemma describes the relation between the systems $\{\varphi_n\}_{n=0}^\infty$ and $\{\varphi_n^l\}_{n=0}^\infty$, $l\in\mathbb N$. 
\begin{lemma}
	\label{lemzeta}
	Let $n\in \mathbb N_0$ and $l>(n+1)/T$. Then 
	\begin{equation}
	\label{normphi}
	\big\|\varphi_n-\varphi_n^l\big\|_{\LL}
	\leq  \frac{n+1}{2^{n+5/2}l} \frac{\sqrt{(2n+2)!}}{(T-(n+1)/l)^{n+3/2}}
	\underset{l\rightarrow\infty}\rightarrow 0.
	\end{equation}
\end{lemma}
\begin{proof}
	Let  $l>\frac{n+1}T$, $n\in \mathbb N_0$. Since
	\begin{align*}
	(1+y)^{n+1}-1&\leq (n+1)y(1+y)^n,& y&\in\mathbb R_+,
	\\
	e^\xi-1-\xi&\leq \frac12\xi^2 e^\xi\quad \text{and}\quad	e^\xi-1\leq \xi e^\xi, & \xi& \in\mathbb R_+,
	\end{align*}
	we have 
	\begin{align*}
	\big|\varphi_n(\rho)-\varphi_n^l(\rho) \big|
	&=\left(\left(\frac{e^{\rho/l}-1}{\rho/l}\right)^{n+1}-1\right) \rho^n e^{-T\rho}
	\\
	&\leq (n+1)\left(\frac{e^{\rho/l}-1}{\rho/l}\right)^n \left(\frac{e^{\rho/l}-1}{\rho/l}-1\right)\rho^n e^{-T\rho}
	\\
	&\leq \frac{n+1}{2l} \rho^{n+1} e^{-(T-(n+1)/l)\rho}
	,\quad \rho \in\mathbb R_+.
	\end{align*}
	Therefore,	
	\begin{align*}
	\big\|\varphi_n-\varphi_n^l\big\|_{\LL}
	&\leq \frac{n+1}{2l} \left(\int_0^\infty \rho^{2(n+1)} e^{-2(T-(n+1)/l)\rho}\,d\rho \right)^{1/2} 
	\\
	&\leq \frac{n+1}{2^{n+5/2}l} \frac{\sqrt{(2n+2)!}}{(T-(n+1)/l)^{n+3/2}}
	\underset{l\rightarrow\infty}\rightarrow 0 ,\quad n\in \mathbb N_0,
	\end{align*}
	that was to be proved.	
\end{proof}


\begin{theorem}
	\label{pthclos}
	Let $T>0$. We have $\overline{\RR}=\LL$.
\end{theorem}

To prove the theorem, we need to construct controls $\{u_n\}_{n=0}^\infty$ from $L^\infty[0,T]$ such that 
\begin{equation}
\EuScript Y_{u_n}(\cdot,T)\underset{n\to\infty}\to g\quad \text{in}\ \LL
\end{equation}
for a given $g\in\LL$.

To this aid, we need an appropriate basis in $\LL$. 
Consider the Laguerre polynomials \cite[pp.~773--775, 22.1.1, 22.1.2, 22.2.13]{HB}:
\begin{equation}
\label{laguer}
L_n(x)
=\frac{e^x}{n!}\left(\frac{d}{dx}\right)^n(e^{-x}x^n)
=\sum_{k=0}^{n}\binom{n}{k}\frac{(-1)^k}{k!}x^k,\quad n\in \mathbb N_0.
\end{equation}
It is well-known that the system $\{\mathrm e_n\}_{n=0}^\infty$ is an orthonormal basis in $\LL$ where $\mathrm e_n(x)=L_n(x)e^{-x/2}$, $x\in\mathbb R_+$. Put
\begin{align}
\label{q}
\psi_n(r)&=\frac1{\sqrt{2T}}L_n\left(\frac r{2T}\right)e^{-\frac{r}{4T}},& r&\in\mathbb R_+,& n&\in\mathbb N_0,
\\
\label{hq}
\widehat \psi_n(\rho)&=(-1)^n \sqrt{2T} L_n(2T\rho)e^{-T\rho},& \rho&\in\mathbb R_+,&n&\in\mathbb N_0.
\end{align}
Evidently, each of the systems $\{\psi_n\}_{n=0}^\infty$ and $\{\widehat \psi_n\}_{n=0}^\infty$ is an orthonormal basis in $\LL$. In addition, we have
\begin{equation}
\label{basis}
\widehat \psi_n=\Phi  \psi_n,\quad n\in\mathbb N_0.
\end{equation}
Let us prove this formula. Let $n\in\mathbb N_0$. With regard to Theorem \ref{prop-2-ii} and \eqref{int-e}, we obtain
\begin{align}
\label{phi-q}
\big(\Phi  \psi_n\big)(\rho)
&=\frac1{2\sqrt{2T}}\int_0^\infty J_0\big(\sqrt{r\rho}\big) L_n\left(\frac r{2T}\right)e^{-\frac{r}{4T}}\,  dr
\notag\\
&=\frac1{2\sqrt{2T}} \sum_{k=0}^n\binom n k \frac{(-1)^k}{k!(2T)^k} \sum_{m=0}^\infty \frac{(-1)^m \rho^m}{(m!)^2 2^{2m}}\int_0^\infty r^{m+k} e^{-\frac{r}{4T}}\,  dr
\notag\\
&=\sqrt{2T}\sum_{k=0}^n\binom n k \frac{(-1)^k 2^k}{k!} \sum_{m=0}^\infty \frac{(-1)^m (T\rho)^m}{m!} \frac{(m+k)!}{m!},\ \rho \in \mathbb R_+.
\end{align}
Due to \eqref{laguer}, we have for the series with $\xi=T\rho$ that
\begin{align*}
\sum_{m=0}^\infty \frac{(-1)^m }{m!} \frac{(m+k)!}{m!} \xi^m
&=\left(\frac{d}{d\xi}\right)^k\sum_{m=0}^\infty \frac{(-1)^m }{m!} \xi^{m+k}
\\
&=\left(\frac{d}{d\xi}\right)^k \big(\xi^k e^{-\xi}\big)
=k! L_k(\xi) e^{-\xi},\quad \xi\in\mathbb R.
\end{align*}
Taking this into account and continuing \eqref{phi-q}, we obtain
$$
\big(\Phi  \psi_n\big)(\rho)
=\sqrt{2T}\sum_{k=0}^n\binom n k (-1)^k 2^k L_k(T\rho) e^{-T\rho},\quad \rho \in \mathbb R_+.
$$
By using multiple argument formula (see \cite[p. 785, 22.12.7]{HB}): 
$$
L_n(\mu \xi)=\sum_{k=0}^{n}\binom{n}{k} \mu^k (1-\mu)^{n-k} L_k(\xi), \quad \mu\in \mathbb R,\ \xi \in \mathbb R,
$$
with $\mu=2$ and $\xi=T\rho$, we get \eqref{hq}.

\begin{proof}[Proof of Theorem \textup{\ref{pthclos}}]
	Let $g\in\LL$. Put $G=\Phi g$. Hence, $G\in\LL$. 
	Set $g_n=\langle g,\psi_n \rangle_{\LL}$, $n\in\mathbb N_0$. We have
	$$
	g=\sum_{n=0}^\infty g_n \psi_n
	$$
	and
	\begin{equation*}
	G=\sum_{n=0}^\infty g_n \widehat \psi_n 
	=G^N + \sum_{n=N+1}^\infty g_n \widehat \psi_n,
	\end{equation*}
	where $N\in\mathbb N$ and
	\begin{equation*}
	G^N= \sum_{n=0}^N g_n \widehat \psi_n.
	\end{equation*}
	Then
	\begin{equation}
	\label{exp1}
	\big\|G-G^N \big\|_{\LL} =\left(\sum_{n=N+1}^\infty|g_n|^2\right)^{1/2}
	\underset{N\to\infty}\to 0.
	\end{equation}
	With regard to \eqref{hq}, we get
	\begin{equation*}
	G^N
	=\sqrt{2T} \sum_{n=0}^N g_n (-1)^n 
	\sum_{k=0}^n \binom nk \frac{(-1)^k}{k!}(2T)^k \varphi_k
	=\sqrt{2T} \sum_{k=0}^N  \frac{(-1)^k}{k!}(2T)^k d_k^N \varphi_k,
	\end{equation*}
	where
	\begin{equation}
	\label{dd}
	d_k^N=\sum_{n=k}^N \binom nk (-1)^n g_n,\quad k=\overline{0,N}.
	\end{equation}
	Put
	\begin{equation}
	\label{exp3}
	G^N_l=\sqrt{2T} \sum_{k=0}^N  \frac{(-1)^k}{k!}(2T)^k d_k^N \varphi_k^l, \quad l\in\mathbb N. 
	\end{equation}
	Taking into account Lemma \ref{lemzeta},  we conclude that  for $l> (N+1)/T$, we have
	\begin{align}
	\label{exp2}
	\big\| G^N-G^N_l\big\|_{\LL}
	&\leq \sqrt{2T} \sum_{k=0}^N  \frac{(2T)^k}{k!} \big|d_k^N\big| \big\|\varphi_k^l-\varphi_k\big\|_{\LL}
	\notag\\
	&\leq \frac{\sqrt T}{4l} \sum_{k=0}^N \frac{T^k \sqrt{(2k+2)!}}{(T-(k+1)/l)^{k+3/2}} \frac{k+1}{k!} \big|d_k^N \big|\underset{l\to\infty}\to 0.
	\end{align}

	Let $g^N_l=\Phi^{-1} G^N_l$. We have	
	\begin{align}
	\label{exp22}
	\left\| g-g^N_l \right\|_{\LL}
	&=\left\| G-G^N_l \right\|_{\LL}
	\notag\\
	&\leq \left\|G-G^N \right\|_{\LL}+\left\| G^N-G^N_l \right\|_{\LL}.
	\end{align}
	With regard to \eqref{exp1} and \eqref{exp2},  we conclude that for all $\varepsilon>0$, we can choose  appropriate $N\in\mathbb N$ and $l> (N+1)/T$ such that 
	\begin{equation}
	\label{exp4}
	\left\| g-g^N_l \right\|_{\LL}<\varepsilon.
	\end{equation} 
	Let us prove that $g^N_l\in\RR$.
	Put 
	\begin{equation}
	u_l^n(\xi)=\begin{cases}
	(-1)^{n-j}\binom{n}{j} l^{n+1},& \xi\in\left(\frac{j}{l},\frac{j+1}{l}\right),\
	j=\overline{0,n}\\
	0,& \xi\notin\left[ 0,\frac{n+1}{l}\right]
	\end{cases},\quad l\in\mathbb{N},\ n\in \mathbb N_0.
	\label{contrl}
	\end{equation}
	Note that $u_l^n\underset{l\to\infty}\to(-1)^n\delta^{(n)}$ in $H^{-1}(\mathbb{R})$ for  each $n\in \mathbb N_0$. Taking into account \eqref{moment6} and \eqref{contrl}, it is easy to obtain $\Phi\EuScript Y_{u_l^n}(\cdot,T)=-\frac{1}{\pi}\varphi_n^l$. Due to \eqref{exp3}, we have 
	\begin{align}
	\label{func}
	g^N_l
	&=-\sqrt{2T}\pi\sum_{k=0}^N  \frac{(-1)^k}{k!}(2T)^k d_k^N \EuScript{Y}_{u_l^k}(\cdot,T) = \EuScript{Y}_{\EuScript{U}^N_l}(\cdot,T),
	\end{align}
	where 
	\begin{equation}
	\label{contr2}
	\EuScript{U}^N_l(\xi)=-\sqrt{2T}\pi\sum_{k=0}^N  \frac{(-1)^k}{k!}(2T)^k d_k^N u_l^k(\xi),\quad \xi\in\mathbb{R}_+.
	\end{equation}
	In addition, due to \eqref{exp1}, \eqref{exp2}, and \eqref{exp22}, we obtain
	\begin{align}
	\label{est}
	&\left\| g-g^N_l \right\|_{\LL}=\left\| g-\EuScript{Y}_{\EuScript{U}^N_l}(\cdot,T) \right\|_{\LL}
	\leq \left(\sum_{n=N+1}^\infty|g_n|^2\right)^{1/2}
	\notag\\
	&\quad
	+ \frac{\sqrt T}{4l} \sum_{k=0}^N \frac{T^k \sqrt{(2k+2)!}}{(T-(k+1)/l)^{k+3/2}} \frac{k+1}{k!} \big|d_k^N \big|,\quad N\in\mathbb N,\ l>\frac{N+1}T.
	\end{align}
	Evidently, $\EuScript{U}^N_l\in L^\infty(0,T)$. Thus, with regard to \eqref{re1} and \eqref{func}, we can see that $g^N_l\in\RR$.
	Since we have considered an arbitrary $\varepsilon>0$ in \eqref{exp4}, we conclude that $g\in\overline{\RR}$.
\end{proof} 





\section{Example}
\label{ex}
In this section, we give an example illustrating Theorem \ref{thappc}.
\begin{example}
	\label{ex111}
	Let 
	$$
	w^0(x)=\cosh\left(\frac{|x|^2}{12T}\right)e^{-\frac{|x|^2}{4T}},\quad w^T(x)=\frac{3}{14}e^{-\frac{|x|^2}{7T}},\quad x\in \mathbb{R}^2.
	$$
	Consider the problem of approximate controllability for system \eqref{eq1}, \eqref{ic1} with $W^0=w^0$ and $W^T=w^T$. 
	
	We have $W^0\in\widehat{H}^0\left(\mathbb R^2\right)$ and $W^T\in\widehat{H}^0\left(\mathbb R^2\right)$. Moreover, $W^0\in\EuScript{H}$ and $W^T\in\EuScript{H}$.
	With regard to \eqref{sol1w0} and \eqref{www}, it is easy to see that
	\begin{align*}  
	\EuScript W_0(x,T)&=\frac{3}{10}e^{-\frac{|x|^2}{10T}}+\frac{3}{14}e^{-\frac{|x|^2}{7T}},&& x\in\mathbb{R}^2,\\
	W^T_0(x,T)&=-\frac{3}{10}e^{-\frac{|x|^2}{10T}},&& x\in\mathbb{R}^2.
	\end{align*}
	Evidently, $W^T_0\in\EuScript{H}$. 
	
	According to Theorem \ref{thappc}, the initial state $W^0$ is approximately controllable to the target state $W^T$ in the given time $T$. Note that condition \eqref{nec} does not hold for $W^T_0$.
	
	To construct controls solving the approximate controllability problem for control system \eqref{eq1}, \eqref{ic1}, we use the method given in the proof of Theorem \ref{pthclos}.
	Put $g=\Psi W_0^T$. Then 
	$$
	g(r)=-\frac3{10} e^{-\frac{r}{10T}},\quad r\in\mathbb R_+.
	$$ 
	With regard to \eqref{int-e}, we get 
	\begin{align*}
	g_n=\langle g, \psi_n \rangle_{\LL}
	&=-\frac{3}{10\sqrt{2T}} \sum_{k=0}^n \binom nk \frac{(-1)^k}{k! (2T)^k} \int_0^\infty r^k e^{-\frac{7r}{20T}}\, dr
	\\
	&= -\frac{3}{10\sqrt{2T}} \sum_{k=0}^n \binom nk \frac{(-1)^k}{k! (2T)^k} k! \left(\frac{20T}{7}\right)^{k+1}
	\\
	&=-\frac37\sqrt{2T} \sum_{k=0}^n \binom nk \left(-\frac{10}{7}\right)^k
	\\
	&=(-1)^{n+1} \left(\frac37\right)^{n+1}\sqrt{2T},\quad  n\in \mathbb N_0. 
	\end{align*}
	
	Therefore,
	\begin{align}
	\label{infty}
	&\left(\sum_{n=N+1}^\infty|g_n|^2\right)^{1/2}
	=\sqrt{2T} \left(\sum_{n=N+1}^\infty\left(\frac{3}{7}\right)^{2n+2}\right)^{1/2}
	\notag\\
	&\qquad
	=\sqrt{2T} \left(\frac{3}{7}\right)^{N+2} \left(\frac{1}{1-9/49}\right)^{1/2}
	= \frac32\sqrt{\frac{T}{5}}\left(\frac{3}{7}\right)^{N+1}, \quad N\in\mathbb N.
	\end{align}
	Due to  \eqref{dd}, we have
	\begin{equation*}
	d_k^N=\sum_{n=k}^N \binom nk (-1)^n g_n
	=-\frac37\sqrt{2T}\sum_{n=k}^N \binom nk  \left(\frac37\right)^n
	,\quad k=\overline{0,N}.
	\end{equation*}
	%
	%
	%
	%
	%
	%
	%
	\begin{figure}[!h]
		\begin{center}
			\begin{subfigure}[t]{0.45\linewidth}
				\centering
				\includegraphics{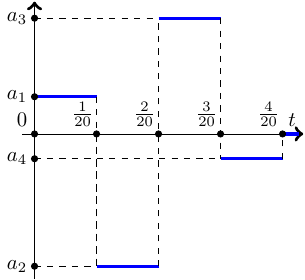}\\
				\centering \caption{\parbox[t]{0.7\textwidth}{ 
						$N=3$, $l=20$,\\ 
						$a_1\approx 4171487.587754723$,\\ 
						$a_2\approx -11985246.36814925$,\\ 
						$a_3\approx 11476859.47814512$,\\ 
						$a_4\approx -3662827.493025041$.
				}}
				\label{fig:u-3-1}
			\end{subfigure}
			\begin{subfigure}[t]{0.53\linewidth}
				\centering
				\includegraphics{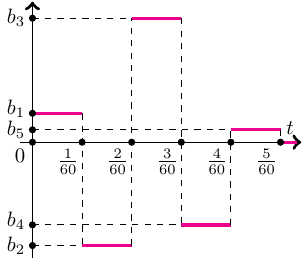}\\
				\centering \caption{\parbox[t]{0.7\textwidth}{
						$N=4$, $l=60$, \\
						$b_1\approx 12268766670.45946$,\\ 
						$b_2\approx -48230066041.31739$, \\
						$b_3\approx 71097757825.27233$, \\
						$b_4\approx -46580177228.79937$,\\ 
						$b_5\approx 11443719610.35109$. 
				}}
				\label{fig:u-3-2}
			\end{subfigure}
			\centering \caption{The controls $\EuScript{U}^N_l$ defined by \eqref{contr3}.} \label{fig:appr}
		\end{center}
	\end{figure}%
	%
	\begin{figure}[!h]
		\begin{center}
			\begin{subfigure}[t]{0.48\linewidth}
				\centering \includegraphics[height=54mm]{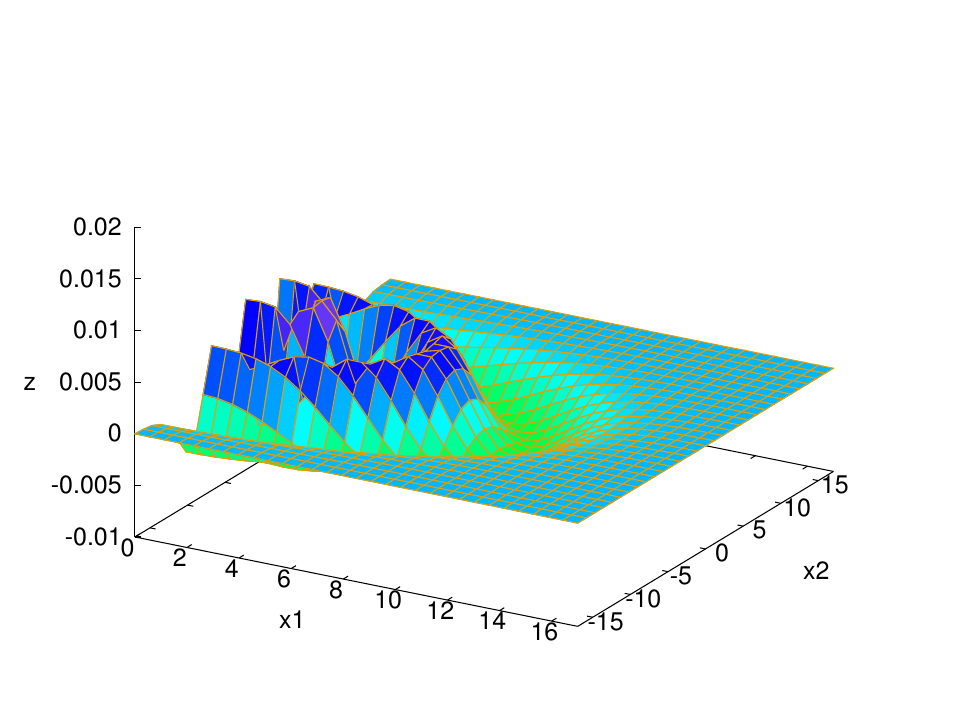}\\
				\centering \caption{\parbox[t]{0.88\textwidth}{$N=3$, $l=20$.}}
				\label{fig:appr11}
			\end{subfigure}
			\begin{subfigure}[t]{0.48\linewidth}
				\centering \includegraphics[height=54mm]{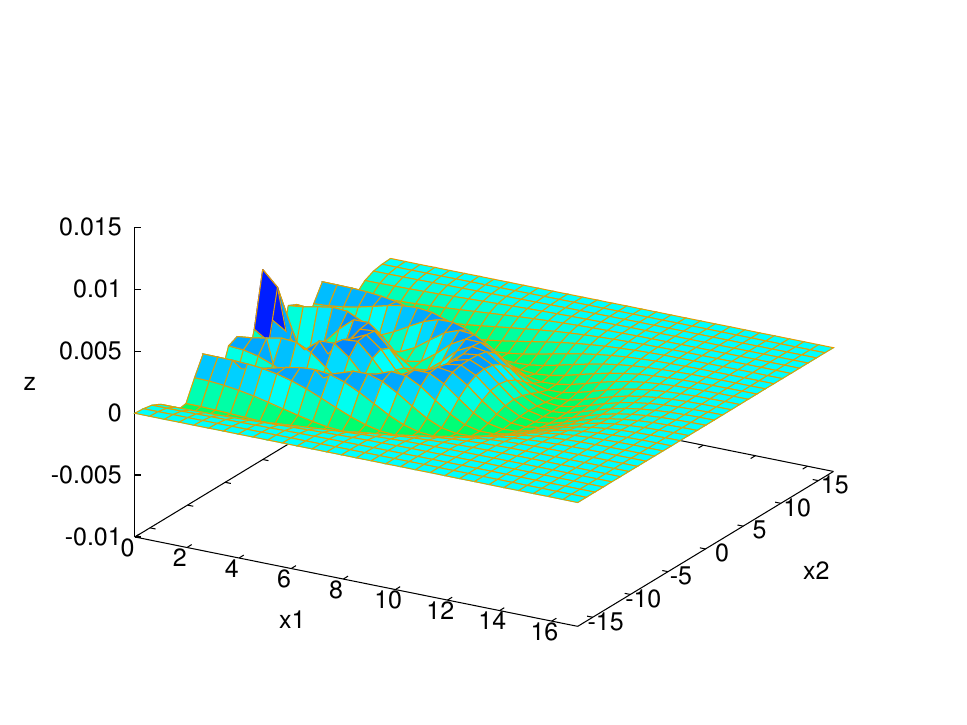}
				\centering \caption{\parbox[t]{0.88\textwidth}{%
						$N=4$, $l=60$.}}
				\label{fig:appr12}
			\end{subfigure}
			\centering \caption{The influence of the controls $\EuScript{U}^N_l$
				on the difference $W^T-\left(\EuScript W_0(\cdot,T)+\EuScript W_{\EuScript{U}^N_l}(\cdot,T)\right)$ with $T=3$.}
			\label{fig:appr1}
		\end{center}
	\end{figure}%
	%
	\begin{figure}[!h]
		\begin{center}
			\begin{subfigure}[t]{0.49\linewidth}
				\centering \includegraphics[height=47mm]{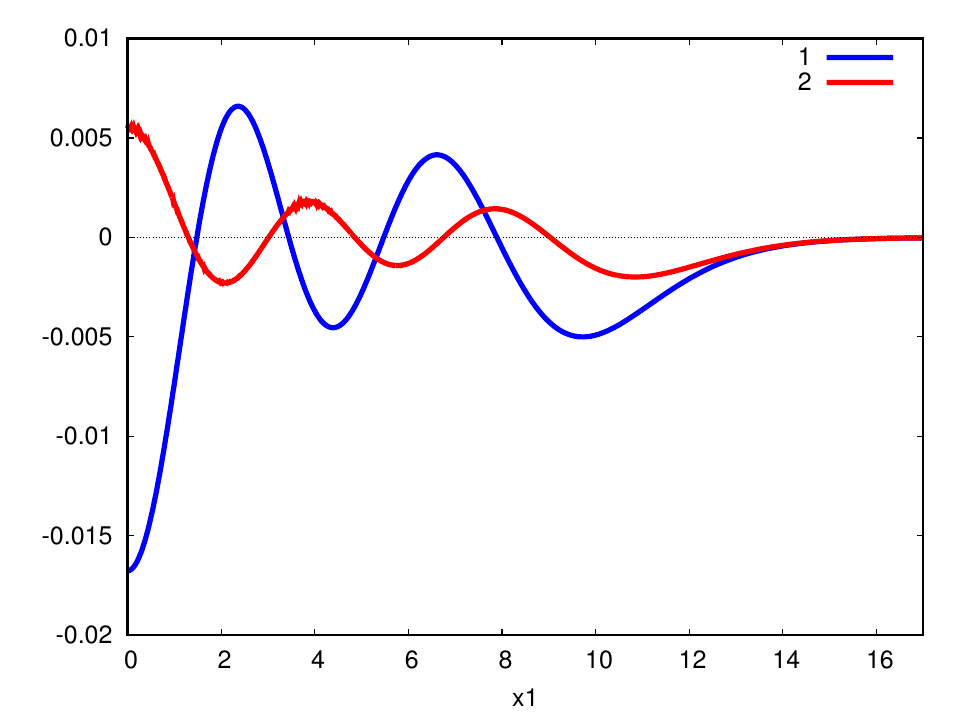}\\
				\centering \caption{\parbox[t]{0.88\textwidth}{
						1) $N=3$, $l=20$; 
						2) $N=4$, $l=60$.}}
				\label{fig:appr21}
			\end{subfigure}
			\begin{subfigure}[t]{0.49\linewidth}
				\centering \includegraphics[height=47mm]{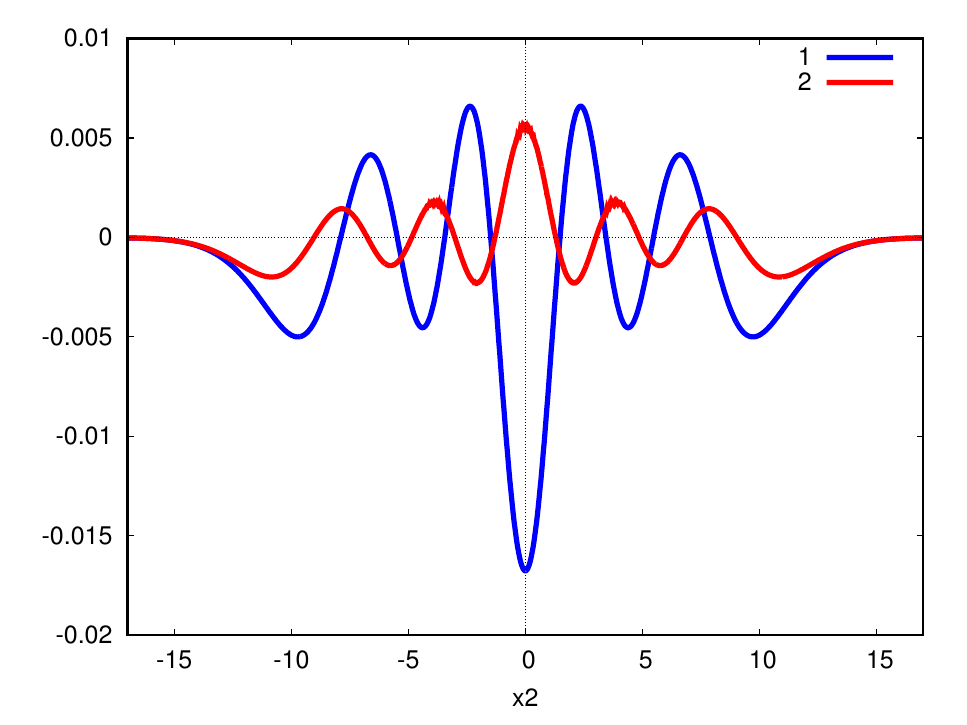}
				\centering \caption{\parbox[t]{0.88\textwidth}{%
						1) $N=3$, $l=20$; 
						2) $N=4$, $l=60$.}}
				\label{fig:appr22}
			\end{subfigure}
			\centering \caption{The influence of the controls $\EuScript{U}^N_l$
				on the difference $W^T-\left(\EuScript W_0(\cdot,T)+\EuScript W_{\EuScript{U}^N_l}(\cdot,T)\right)$  with $T=3$ (vertical section for $x_2=0$ and horizontal section for $x_1=0$).}
			\label{fig:appr2}
		\end{center}
	\end{figure}%
	It follows from \eqref{contr2} that
	\begin{align}
	\label{contr3}
	\EuScript{U}^N_l&=-\sqrt{2T}\pi\sum_{k=0}^N  \frac{(-1)^k}{k!}(2T)^k d_k^N u_l^k
	\notag\\
	&=2T\pi\sum_{k=0}^N \sum_{n=k}^N \binom nk \left(\frac{3}{7}\right)^{n+1} \frac{(-1)^k}{k!}(2T)^k u_l^k.
	\end{align}
	Taking into account \eqref{est} and \eqref{infty}, we obain
	\begin{align*}
	&\big\| g-\EuScript{Y}_{\EuScript{U}^N_l}(\cdot,T) \big\|_{\LL}
	\leq \frac 32 \sqrt{\frac{T}{5}}\left(\frac{3}{7}\right)^{N+1}
	\\
	&+\frac{T}{2\sqrt 2 l} \sum_{k=0}^N \frac{T^k \sqrt{(2k+2)!}}{(T-(k+1)/l)^{k+3/2}} \frac{k+1}{k!} 
	\sum_{n=k}^N \binom nk \left(\frac{3}{7}\right)^{n+1},\  N\in\mathbb N,\ l>\frac{N+1}T.
	\end{align*}
	We have $\EuScript{W}_{\EuScript{U}^N_l}(\cdot,T)=\Psi^{-1}\EuScript{Y}_{\EuScript{U}^N_l}(\cdot,T)$ and $W_0^T=\Psi^{-1} g$. Taking into account Proposition \ref{prop-1-iv}, we get
	\begin{align*}
	\big\|W^T-\big(\EuScript W_0 (\cdot,T)&+\EuScript W_{\EuScript{U}^N_l}(\cdot,T)\big)\big\|^0
	=\|W_0^T-\EuScript{W}_{\EuScript{U}^N_l}(\cdot,T)\|^0
	\\
	&=\sqrt\pi \big\| g-\EuScript{Y}_{\EuScript{U}^N_l}(\cdot,T) \big\|_{\LL},\quad  N\in\mathbb N,\ l>\frac{N+1}T.
	\end{align*}

	The  plots of the controls $\EuScript{U}^N_l$ are given in Fig. \ref{fig:appr} for $T=3$ with the cases of $N=3$, $l=20$ and $N=4$, $l=60$.  
	Figs. \ref{fig:appr1} and \ref{fig:appr2} describe the influence of the control $\EuScript{U}^N_l$ on the difference $W^T-\left(\EuScript W_0(\cdot,T)+\EuScript W_{\EuScript{U}^N_l}(\cdot,T)\right)$ with  $T=3$. 
\end{example}

\section{Proofs of theorems and auxiliary statements}
\label{sect6}

\hypertarget{th34}{}
\begin{proof}[Proof of Theorem \textup{\ref{th-sol}}]
	We prove the theorem using similar reasoning to those of the corresponding theorem in \cite{FKh3}.
	Put $V^0=\mathscr FW^0$, $\EuScript V_0(\cdot,t)=\mathscr F_{x\to\sigma} \EuScript W_0(\cdot,t)$, $\EuScript V_u(\cdot,t)=\mathscr F_{x\to\sigma} \EuScript W_u(\cdot,t)$, $t\in[0,T]$. Evidently,
	\begin{align}
	\label{a2}
	\EuScript V_0(\sigma,t) &= e^{-t|\sigma|^2} V^0(\sigma),&&\sigma\in\mathbb R^2,\ t\in[0,T],
	\\
	\label{a3}
	\EuScript V_u(\sigma,t)&=-\frac{1}{\pi} 
	\int_0^t e^{-\xi|\sigma|^2} u(t-\xi)\, d\xi,&&\sigma\in\mathbb R^2,\ t\in[0,T].
	\end{align}
	Hence,
	\begin{equation}
	\label{a4}
	\|\EuScript W_0(\cdot,t)\|^0=\|\EuScript V_0(\cdot,t)\|_0\leq \|V^0\|_0=\|W^0\|^0,\quad t\in[0,T].
	\end{equation}
	Thus, \ref{th-sol-ia} is proved.
	
	Let $\alpha=(\alpha_1,\alpha_2)\in\mathbb N^2_0$. Then
	\begin{align}
	\label{a5}
	\left\| D^\alpha \EuScript W_0(\cdot,t) \right\|^0
	&=\left\|(\cdot)_{[1]}^{\alpha_1}(\cdot)_{[2]}^{\alpha_2} \EuScript V_0(\cdot,t) \right\|_0
	\leq e^{ t} \left(\frac{1+|\alpha|}{2t e}\right)^{(1+|\alpha|)/2} \|V^0\|_{-1}
	\nonumber\\
	&\leq e^{t} \left(\frac{1+|\alpha|}{2t e}\right)^{(1+|\alpha|)/2} \|W^0\|^0
	,\quad  t\in(0,T].
	\end{align}
	Here we have used the following estimates:
	\begin{align*}
	&\left|\sigma_1^{\alpha_1}\sigma_2^{\alpha_2}\EuScript V_0(\sigma,t) \right|^2\leq \left(1+|\sigma|^2\right)^{1+|\alpha|} e^{-2t|\sigma|^2} \frac{|V^0(\sigma)|^2}{1+|\sigma|^2},&& \sigma \in \mathbb R^2,\ t\in[0,T],\\
	&\xi^m e^{-\beta \xi}\leq \left(\frac{m}{\beta e}\right)^m,&& m\in\mathbb N,\  \beta>0,\ \xi\geq 0,
	\\
	&\|V^0\|_{-1}\leq \|V^0\|_0=\|W^0\|^0.
	\end{align*}
	Thus, \ref{th-sol-iia} is proved.
	
	Suppose $W^0\in \EuScript H$. Since $\EuScript W_0(\cdot,t)=\mathscr F^{-1}_{\sigma\to x} \EuScript V_0(\cdot,t)$, with regard to \eqref{a2}, we get
	$\EuScript W_0(\cdot,t)\in \EuScript H$, $t\in[0,T]$. 
	Thus, \ref{th-sol-iiia}  holds.
	
	It remains to prove \ref{th-sol-iva}. Put
	$$
	g(r,t)=\int_0^t e^{-\xi r^2} u(t-\xi)\, d\xi,\quad r\geq 0, \ t\in[0,T].
	$$
	Taking into account, that 
	$$
	\frac{1-e^{-tr^2}}{r^2}\leq \frac{2(t+1)}{r^2+1},\quad r>0,\ t>0,
	$$
	we obtain
	\begin{align*}
	|g(r,t)|&\leq \|u\|_{L^\infty(0,T)} \frac{1-e^{-tr^2}}{r^2}
	\leq  \|u\|_{L^\infty(0,T)} \frac{2(t+1)}{r^2+1},\quad r>0,\ t>0.
	\end{align*}
	Hence,
	\begin{align*}
	&\left\| \EuScript W_u(\cdot,t) \right\|^0
	= \left\| \EuScript V_u(\cdot,t) \right\|_0
	=\sqrt{\frac{2}{\pi}}\left(\int_0^\infty |g(r,t)|^2r\, dr \right)^{1/2}\\
	&\leq \frac{2}{\sqrt{\pi}}(t+1)  \|u\|_{L^\infty(0,T)} \left(\int_0^\infty \frac{2r\, dr}{(1+r^2)^2}\right)^{1/2}
	=\frac{2}{\sqrt{\pi}}(t+1)  \|u\|_{L^\infty(0,T)},\  t\in[0,T].
	\end{align*} 
	This completes the proof.
\end{proof}

\begin{lemma}
	\label{lempoint}
	Let $W^0\in \widehat{H}^0\left(\mathbb{R}^2\right)$, $t\in [0,T]$. Let $W$ be a solution to \eqref{eq1}, \eqref{ic1}. Then \eqref{bc1} holds.
\end{lemma}

\begin{proof}
	We have from \eqref{sol11}
	\begin{equation}
	\label{sol1dif}
	\frac{\partial}{\partial x_1}W(0^+,(\cdot)_{[2]},t)=\frac{\partial}{\partial x_1} \EuScript W_0(0^+,(\cdot)_{[2]},t)+\frac{\partial}{\partial x_1}\EuScript W_u(0^+,(\cdot)_{[2]},t),\ t\in(0,T].
	\end{equation}
	According to Theorem 
	\ref{th-sol-ii}, $\frac{\partial}{\partial x_1}\EuScript W_0(\cdot,t)$ is  continuous on $\mathbb R^2$ for each $t\in(0,T]$. Moreover, $\frac{\partial}{\partial x_1}\EuScript W_0(\cdot,t)$ is odd with respect to $x_1$, $t\in[0,T]$. Hence,
	\begin{equation}
	\label{aa1}
	\frac{\partial}{\partial x_1}\EuScript W_0(0^+,(\cdot)_{[2]},t)=0,\quad t\in(0,T].
	\end{equation}
	For $\frac{\partial}{\partial x_1}\EuScript W_u(0^+,(\cdot)_{[2]},t)$, $t\in(0,T]$ we have
	$$
	\frac{\partial}{\partial x_1}\EuScript W_u(x,t)=\frac2\pi \frac{x_1}{|x|^2} \int_{\frac{|x|}{2\sqrt t}}^\infty y e^{-y^2} 
	u\left(t-\frac{|x|^2}{4y^2}\right)\, dy, \quad x\in\mathbb R^2,\ t\in(0,T].
	$$
	Let $\varphi\in H^1(\mathbb R)$, $\psi\in\EuScript D(\mathbb R_+)$, and $\alpha>0$. 
	With regard to Definition \ref{defpoint}, we consider \eqref{value} for $f=\frac{\partial}{\partial x_1}\EuScript W_u(\cdot,t)$, $t\in(0,T]$, changing  additionally the variable $x_j$ to the new variable $\alpha x_j$ in the integral with respect to $x_j$, $j=1,2$.
	We have
	\begin{align}
	\label{aa2}
	&\left\langle\left\langle \frac{\partial}{\partial x_1}\EuScript W_u( (\cdot)_{[1]},(\cdot)_{[2]}),\varphi((\cdot)_{[2]})\right\rangle_{[2]},\frac1\alpha\psi\left(\frac{(\cdot)_{[1]}}\alpha\right)\right\rangle_{[1]}
	\nonumber\\
	&\quad= \frac2\pi \int_{-\infty}^\infty\int_{-\infty}^\infty \frac{x_1}{|x|^2} 
	\int_{\frac{\alpha|x|}{2\sqrt t}}^\infty y e^{-y^2} 
	u\left(t-\frac{\alpha^2|x|^2}{4y^2}\right) dy\, \varphi(\alpha x_2)\, dx_2\, \psi(x_1)\, dx_1,
	\nonumber\\
	&\kern60ex t\in(0,T].
	\end{align}
	Since
	\begin{align*}
	\int_{\frac{\alpha|x|}{2\sqrt t}}^\infty y e^{-y^2} 
	\left|u\left(t-\frac{\alpha^2|x|^2}{4y^2}\right)\right| dy
	\leq \|u\|_{L^\infty(0,T)} \int_0^\infty y e^{-y^2}\, dy = \frac{1}{2}\|u\|_{L^\infty(0,T)}
	\end{align*}
	and
	\begin{align*} 
	\int_{-\infty}^\infty\int_{-\infty}^\infty \frac{x_1}{|x|^2} |\varphi(\alpha x_2)|&\, dx_2\, |\psi(x_1)|\, dx_1
	\\
	&
	\leq \sup_{\mu\in\mathbb R} |\varphi(\mu)| \int_0^\infty x_1 |\psi(x_1)| \int_{-\infty}^\infty \frac{dx_2}{x_1^2+x_2^2}\,dx_1
	\\
	&=\pi  \sup_{\mu\in\mathbb R} |\varphi(\mu)| \int_0^\infty  |\psi(x_1)|\, dx_1<\infty,
	\end{align*}
	we can apply Lebesgue's dominated convergence theorem to \eqref{aa2} as $\alpha\to0^+$:
	\begin{align*}
	&\left\langle\left\langle \frac{\partial}{\partial x_1}\EuScript W_u( (\cdot)_{[1]},(\cdot)_{[2]}),\varphi((\cdot)_{[2]})\right\rangle_{[2]},\frac1\alpha\psi\left(\frac{(\cdot)_{[1]}}\alpha\right)\right\rangle_{[1]}
	\\
	&\qquad\to \frac2\pi u(t) \varphi(0)
	\int_{0}^\infty\int_{-\infty}^\infty\frac{x_1}{|x|^2} 
	\int_{0}^\infty y e^{-y^2} 
	\, dy\, dx_2\, \psi(x_1)\, dx_1
	\\
	&\qquad=  \frac1\pi u(t) \varphi(0)  \int_0^\infty x_1 \psi(x_1) \int_{-\infty}^\infty \frac{dx_2}{x_1^2+x_2^2}\,dx_1= u(t) \varphi(0)  \int_{-\infty}^\infty \psi(x_1)\, dx_1
	\\
	&\qquad=\left\langle \left\langle u(t)\delta_{[2]},\varphi\right\rangle_{[2]}, \psi \right\rangle_{[1]},\quad t\in(0,T],
	\end{align*}
	i.e. 
	\begin{equation}
	\label{aa3}
	\frac{\partial}{\partial x_1}\EuScript W_u\big(0^+,(\cdot)_{[2]},t\big)= u(t) \delta_{[2]},\quad t\in(0,T].
	\end{equation}
	With regard to \eqref{sol1dif}, \eqref{aa1}, and \eqref{aa3}, we conclude that  \eqref{bc1} holds.
\end{proof}

\hypertarget{th36}{}
\begin{proof}[Proof of Theorem \textup{\ref{thnec}}] 
	According to Theorem \ref{th-sol-iv}, we have $f\in\HH$. Therefore, $F=\mathscr F f\in\HH$ (see Theorem \ref{prop-1-vi}). Set $g=\Psi f$ and $G=\Psi F$.  
	It follows from Theorem \ref{pthnec} and \eqref{r} that Theorem \ref{thnec} is true.
\end{proof}

\hypertarget{th37}{}
\begin{proof}[Proof of Theorem \textup{\ref{entire}}] 
	According to Theorem \ref{th-sol-iv}, we have $f\in\HH$ (see Theorem \ref{prop-1-vi}). Therefore, $F=\mathscr F f\in\HH$. Set $g=\Psi f$ and $G=\Psi F$.  
	It follows from Theorem \ref{pentire} and \eqref{r} that Theorem \ref{entire} is true. 
\end{proof}

\hypertarget{th38}{}
\begin{proof}[Proof of Theorem \textup{\ref{cond-f}}] 
	Set $g=\Psi f$ and $G=\Psi F$.  
	It follows from Theorem \ref{pcond-f} that Theorem \ref{cond-f} is true. 
\end{proof}

\mbox{}
\hypertarget{th39}{}
\begin{proof}[Proof of Theorem \textup{\ref{thmom}}] 
	Set $g=\Psi W_0^T$. Let us prove that for $\{\omega_n\}_{n=0}^\infty$ defined by  \eqref{moment} and $\{\gamma_n\}_{n=0}^\infty$ defined by \eqref{pmoment}, we have
	\begin{equation}
	\label{omg}
	\omega_n=\gamma_n,\quad n\in\mathbb N_0.
	\end{equation}
	Let $n\in\mathbb N_0$. Setting $x_1=\sqrt r\cos\phi$ and $x_2=\sqrt r\sin\phi$, $r\in\mathbb R_+$, $\phi\in[0,2\pi)$, we have
	\begin{equation*}
	\omega_n=-\frac12 \frac{n!}{(2n)!} \iint_{\mathbb R^2} x_1^{2n} W_0^T(x)\, dx
	=-\frac14\frac{n!}{(2n)!} \int_0^{2\pi}\cos^{2n}\phi\, d\phi \int_0^\infty r^n g(r)\, dr.
	\end{equation*}
	Since
	\begin{align*}
	\frac14 \int_0^{2\pi}\cos^{2n}\phi\, d\phi 
	&= \int_0^\infty \frac{dy}{(1+y^2)^{n+1}}
	= \frac{2n-1}{2n} \int_0^\infty \frac{dy}{(1+y^2)^n}
	\\
	&=\frac{(2n-1)!!}{(2n)!!} \int_0^\infty \frac{dy}{1+y^2}
	=\frac\pi2 \frac{(2n-1)!!}{(2n)!!}
	=\pi \frac{(2n)!}{2^{2n+1}(n!)^2},
	\end{align*}
	we have
	\begin{equation*}
	\omega_n=-\pi \frac{(2n)!}{2^{2n+1}(n!)^2}\frac{n!}{(2n)!}  \int_0^\infty r^n g(r)\, dr
	=-\frac{\pi}{2^{2n+1} n!} \int_0^\infty r^n g(r)\, dr
	=\gamma_n,
	\end{equation*}
	i.e. \eqref{omg} holds.
	It follows from Theorems \ref{pthmom} and \ref{laa4}, \eqref{r}, and \eqref{omg} that Theorem \ref{thmom} is true. 
\end{proof}

\hypertarget{th311}{}
\begin{proof}[Proof of Theorem \textup{\ref{thmomap}}] 
	Set $g=\Psi W_0^T$. For $\{\omega_n\}_{n=0}^\infty$ defined by  \eqref{moment} and $\{\gamma_n\}_{n=0}^\infty$ defined by \eqref{pmoment},
	it follows from Theorems \ref{laa4} and \ref{pthmomap}, \eqref{r}, \eqref{rrr}, and \eqref{omg} that Theorem \ref{thmomap} is true.	
\end{proof}

\hypertarget{th313}{}
\begin{proof}[Proof of Theorem \textup{\ref{thclos}}] 
	It follows from Theorems \ref{prop-1-i} and \ref{pthclos}, \eqref{r}, \eqref{rrr} that Theorem \ref{thclos} is true.	
\end{proof}


\hypertarget{th315}{}
\begin{proof}[Prof of Theorem \textup{\ref{thnulcontr}}]
	Let a state $W^0\in\widehat{H}^0\left(\mathbb R^2\right)$ be controllable to the state $W^T=0$. Then  there exists a control $u\in L^\infty(0,T)$ such that there exists a unique solution $W$ to system \eqref{eq1}, \eqref{ic1} under this control and $W(\cdot,T)=0$. 
	It follows  from \eqref{sol11} that
	\begin{equation*}
	\big(\mathscr F W^0\big)(\sigma)=\frac1\pi\int_0^T e^{\xi|\sigma|^2}u(\xi)\, d\xi,\quad \sigma\in\mathbb{R}^2.
	\end{equation*}
	Evidently, $\mathscr F W^0\in\HH$.
	Setting $G=\Psi \mathscr F W^0$, we obtain
	\begin{equation}
	\label{ccc}
	G(\rho)=\frac1\pi\int_0^T e^{\xi\rho}u(\xi)\, d\xi,\quad \rho\in\mathbb{R}_+.
	\end{equation}
	Let $T^*>T$. Put
	\begin{align}
	\label{hqq}
	\widehat \psi_n^*(\rho)&=(-1)^n \sqrt{2T^*} L_n(2T^*\rho)e^{-T^*\rho},& \rho&\in\mathbb R_+,&n&\in\mathbb N_0,
	\\
	\label{hqq-1}
	\alpha_n&=\langle G,\widehat \psi_n^*\rangle_{\LL}, &&&n&\in\mathbb N_0,
	\\
	\label{hqq-2}
	\beta_n(\xi)&=\frac1\pi \big\langle e^{\xi (\cdot)},\widehat \psi_n^*\big\rangle_{\LL}, &\xi&\in[0,T],&n&\in\mathbb N_0.
	\end{align}
	Obviously,  the system  $\{\widehat \psi_n^*\}_{n=0}^\infty$ is  an orthonormal basis in $\LL$ (cf. \eqref{hq}).
	Then, due to  \eqref{ccc}, we have
	\begin{equation*}
	\sum_{n=0}^\infty\alpha_n \widehat \psi_n^*
	=\sum_{n=0}^\infty\left(\int_0^T\beta_n(\xi)u(\xi)\, d\xi\right) \widehat \psi_n^*.
	\end{equation*}
	Hence,
	\begin{equation}
	\label{ccc3}
	\int_0^T \beta_n(\xi)u(\xi)\,d\xi=\alpha_n,\quad n\in \mathbb N_0.
	\end{equation}
	Let $n\in \mathbb N_0$ be fixed. Taking into account \eqref{laguer}, we get
	\begin{align}
	\beta_n(\xi)&=\frac{(-1)^n}{\pi} \sqrt{2T^*}
	\sum_{k=0}^n\binom nk \frac{(-1)^k}{k!}(2T^*)^k \int_0^\infty \rho^k e^{-(T^*-\xi)\rho}\, d\rho
	\nonumber\\
	&=\frac{(-1)^n}{\pi} \sqrt{2T^*}
	\sum_{k=0}^n\binom nk \frac{(-1)^k}{k!}(2T^*)^k \frac{k!}{(T^*-\xi)^{k+1}}
	\nonumber\\
	&=\frac{(-1)^n}{\pi} \frac{\sqrt{2T^*}}{T^*-\xi} \sum_{k=0}^n \binom nk \left(-\frac{2T^*}{T^*-\xi}\right)^k\nonumber\\
	&=\frac1\pi\frac{\sqrt{2T^*}}{T^*-\xi}\left(\frac{T^*+\xi}{T^*-\xi}\right)^n,
	\quad \xi\in[0,T].
	\label{ccc3a}
	\end{align}
	According to \eqref{ccc3}, we obtain
	\begin{equation*} 
	\alpha_n=\frac{\sqrt{2T^*}}{\pi}\int_0^T\left(\frac{T^*+\xi}{T^*-\xi}\right)^n\frac{u(\xi)}{T^*-\xi}\,d\xi
	=\frac{\sqrt{2T^*}}{\pi}\int_0^{T_1} e^{n\tau}u\left(T^*\frac{e^\tau-1}{e^\tau+1}\right)\frac{e^\tau\, d\tau}{e^\tau+1},
	\end{equation*}
	where $\frac{T^*+\xi}{T^*-\xi}=e^\tau$,  $T_1=\ln\frac{T^*+T}{T^*-T}$.
	Set 
	$$
	\alpha_n^*=\frac{\pi}{\sqrt{2T^*}}\alpha_n,\quad n\in \mathbb N_0,\quad
	u^*(\tau)=\frac{e^\tau}{e^\tau+1}u\left(T^*\frac{e^\tau-1}{e^\tau+1}\right),\quad \tau\in[0,T_1].
	$$ 
	Then,
	\begin{equation}
	\label{ccc4}
	\int_0^{T_1} e^{n\tau}u^*(\tau)\,d\tau=\alpha_n^*,\quad n\in \mathbb N_0.
	\end{equation}
	With regard to \eqref{hqq-1}, we get
	\begin{align*}
	\left|\alpha_n^*\right|
	\leq\frac{\pi}{\sqrt{2T^*}}
	\| G\|_{\LL},\quad n\in \mathbb N_0.
	\end{align*}
	Therefore, for all $\delta>0$, there exists $C_\delta>0$ such that
	\begin{equation}
	\label{ccc5}
	\left| \alpha_n^*\right|\le C_\delta e^{n\delta},\quad n\in \mathbb N_0.
	\end{equation}
	Obviously,
	\begin{equation}
	\label{ccc6}
	\|u^*\|_{L^2(0,T_1)}=\left(\frac{1}{2T^*}\int_0^T \left| u(\xi)\right|^2\,d\xi\right)^{1/2}
	\leq\left(\frac{T}{2T^*}\right)^{1/2}\left\| u\right\|_{L^\infty(0,T)}.
	\end{equation}
	Taking into account \eqref{ccc4}--\eqref{ccc6}, we conclude that all assertions of \cite[Theorem~3.1, b)]{MZua1} hold. Thus, due to this theorem, $\alpha_n^*=0$,
	$n\in \mathbb N_0$. Hence, $G=0$. Therefore, $W^0=0$.
\end{proof}


\subsection*{Acknowledgments.} The authors would like to thank the anonymous referee for their valuable comments and suggestions.
The authors' research is partially supported by the Akhiezer Foundation.



\EndPaper


\end{document}